\documentclass[journal]{IEEEtran}
\IEEEoverridecommandlockouts
\ifCLASSINFOpdf
\else
\fi

\usepackage{mathrsfs}
\usepackage{color}
\usepackage{amsfonts}
\usepackage{amsthm}
\usepackage{subfigure}
\usepackage{caption}
\usepackage{booktabs}
\usepackage{balance}
\usepackage{graphicx}
\usepackage{float} 
\usepackage{enumerate}
\usepackage{algorithm}
\usepackage{algpseudocode}
\usepackage{bm}
\usepackage{makecell}
\usepackage{epstopdf}
\usepackage{mathtools}
\usepackage{bbm}
\usepackage{dsfont}
\usepackage{pifont}
\usepackage{cite}
\usepackage{balance}
\usepackage[T1]{fontenc}
\usepackage{soul}
\usepackage{bbding} % for \checkmark
\usepackage[normalem]{ulem}
\usepackage{tikz}
\newcommand*\circled[1]{\tikz[baseline=(char.base)]{
            \node[shape=circle,draw,inner sep=0.6pt] (char) {\small #1};}}

\newtheorem{remark}{Remark}

\newtheorem{theorem}{Theorem}
\newtheorem{lemma}{Lemma}

\newtheorem{assumption}{Assumption}

\newtheorem{definition}{Definition}

\begin{document}
\title{\LARGE \bf  Distributionally Robust Probabilistic Prediction\\ for Stochastic Dynamical Systems}
\author{Tao Xu and Jianping He
    \thanks{This work was supported by the National Natural Science Foundation of China under Grant 62373247, 625B2117. \textit{(Corresponding author: Jianping He.)} The authors are with the Department of Automation, Shanghai Jiao Tong University, and Key Laboratory of System Control and Information Processing, Ministry of Education of China, Shanghai 200240, China. E-mail: \{Zerken, jphe\}@sjtu.edu.cn.} }

\maketitle

\begin{abstract}
    Probabilistic prediction of stochastic dynamical systems (SDSs) aims to accurately predict the conditional probability distributions of future states. 
    However, accurate probabilistic predictions tightly hinge on accurate distributional information from a nominal model, which is hardly available in practice. 
    To address this issue, we propose a novel functional-maximin-based distributionally robust probabilistic prediction (DRPP) framework. 
    In this framework, one can design probabilistic predictors that have worst-case performance guarantees over a pre-defined ambiguity set of SDSs.
    Nevertheless, DRPP requires optimizing over the space of probability measures with density functions with respect to the Lebesgue measure, which is generally intractable. We develop a methodology that equivalently transforms the original maximin from function spaces to Euclidean spaces. 
    Although it remains intractable to seek a global optimal solution, two suboptimal solutions are derived. By relaxing the constraints on the ambiguity set, we obtain a suboptimal predictor called Noise-DRPP. Relaxing the constraints on the predictor yields another suboptimal predictor, Eig-DRPP. Moreover, optimality gaps between the proposed predictors and the global optimal predictor are derived.
    Finally, we conduct elaborate numerical simulations to compare the performance of different predictors under different SDSs.
\end{abstract}


\section{Introduction}
\subsection{Background}
A probabilistic prediction is a forecast that assigns probabilities to different possible outcomes, enabling better uncertainty quantification, risk assessment, and decision-making. Due to the ever-increasing demand for both accurate and reliable predictions against uncertainty \cite{sevellecNovelProbabilisticForecast2018}, it has wide applications in fields such as epidemiology \cite{cramer2022evaluation}, climatology \cite{buizzaValueProbabilisticPrediction2008}, and robotics \cite{fisac2018probabilistically}, etc. 
The performance of a probabilistic prediction is evaluated by both sharpness and calibration of the predictive probability density function (pdf). 
The sharpness reflects how concentrated (e.g., low-variance or low-entropy) the predictive pdf is, subject to the intrinsic stochasticity of the underlying system. The calibration concerns the statistical compatibility between the predictive pdf and the realized outcomes \cite{gneitingProbabilisticForecasting2014}. To simultaneously assess both sharpness and calibration, one usually assigns a numerical score to the predictive pdf and the realized value of the prediction target. If one has a further requirement that the expected scoring rule is maximized only when the predictive pdf equals the target's real pdf almost everywhere, a strictly proper scoring rule (see in \cite{gneitingStrictlyProperScoring2007a}) is needed.

The probabilistic prediction of a stochastic dynamical system (SDS) is to predict the conditional pdfs of future states based on some prior information about the system. One's prior information is usually in the form of a nominal model, which is a simplified or idealized representation of the underlying system.
While probabilistic prediction for a linear system with Gaussian noises is straightforward, it becomes computationally expensive once the nominal SDS is nonlinear or the noises are non-Gaussian. Therefore, a lot of research contributes to approximating the future pdf to balance the requirement of prediction precision and computation complexity \cite{landgrafProbabilisticPredictionMethods2023}. 

\subsection{Motivation}
Despite the significant achievements in the current prediction algorithms for SDSs, they are not guaranteed to perform well if the distributional information provided by a nominal model is inaccurate. 
It is ubiquitous for a nominal SDS to possess both inaccurate state evolution functions and incomplete distributional information of noises. For example, the nominal state evolution function may be a local linearization from a nonlinear model \cite{hinrichsenUncertainSystems2005}. Additionally, one's nominal information about the distributions of system noises is usually partial. Most of the time, only the estimated mean and covariance are available in the nominal model, which is far from uniquely determining a pdf \cite{schmudgenMomentProblem2017}. 
A nominal predictor can lead to misleading predictive pdfs even though the uncertainties on SDS are not prominent (as demonstrated by experiments in Sec. VIII). Addressing this issue is crucial for the widespread application of probabilistic prediction in model-uncertain scenarios. In this paper, the aim is to develop a framework to design probabilistic predictors for SDSs with worst-case performance guarantees.

To design probabilistic predictors against distributional uncertainties of the model, we inherit the key notions from the distributionally robust optimization (DRO) \cite{delageDistributionallyRobustOptimization2010}. DRO is a minimax-based mathematical framework for decision-making under uncertainty. It aims to find the decision that performs well across a range of possible pdfs in an ambiguity set. Nevertheless, what we need is to find the optimal probabilistic predictor that performs well over a set of SDSs. In DRPP, we generalize the idea of ambiguity sets to jointly describe one's uncertainties of noises' pdfs and state evolution functions. 

The proposed DRPP framework has direct applicability in control and decision-making systems. Reliable predictive pdfs and confidence regions produced by DRPP can be integrated into downstream modules such as risk-aware model predictive control, robot trajectory planning, and safety assessment in power-system frequency regulation. For instance, understanding the evolution of robust confidence regions enables safe motion planning for manipulators and informs constraint tightening in safety-critical control.

\subsection{Contributions}
Motivated by the aforementioned considerations, we propose a novel functional-maximin-based distributionally robust probabilistic prediction (DRPP) framework. 
In DRPP, the objective is the expected score under the worst-case SDS within an ambiguity set, and the predictor optimizes the worst-case objective over the space of admissible predictive pdfs.
% Compared to DRO, the proposed DRPP is significantly different in three perspectives. i) (\textit{Functional optimization}) One needs to find the optimal pdfs over function spaces rather than optimizing decisions over Euclidean spaces. ii) (\textit{Conditional prediction}) The predictor for an SDS is a map from the state space to the pdf space, rather than a single decision variable without conditioning on some context. iii) (\textit{Trajectory expectation}) The objective is an expectation concerning all possible trajectories and control inputs, rather than an expectation over a single random variable. 
The main contributions of this paper are summarized as follows:
\begin{itemize}
    \item To the best of our knowledge, this is the first work on probabilistic prediction of SDSs that optimizes the worst-case performance over an ambiguity set. We have generalized the classical conic moment-based ambiguity set such that both the uncertainties of state evolution functions and system noises are quantified.
    \item We develop a methodology to greatly reduce the complexity of solving DRPP. First, we use the principle of dynamic programming to derive the Bellman equation. Second, we exploit a necessary optimality condition to equivalently transform the maximin problem from function spaces to Euclidean spaces.
    \item We design two suboptimal DRPP predictors by relaxing the constraints. Noise-DRPP is the optimal predictor when the nominal state evolution function is accurate, and Eig-DRPP is the optimal predictor when the predictor is constrained by an eigenvector restriction. Moreover, an upper bound of the optimality gap of the proposed predictors is obtained.
\end{itemize}

\subsection{Organization}
\begin{figure}[htb]
    \centering
    \includegraphics[width=\linewidth]{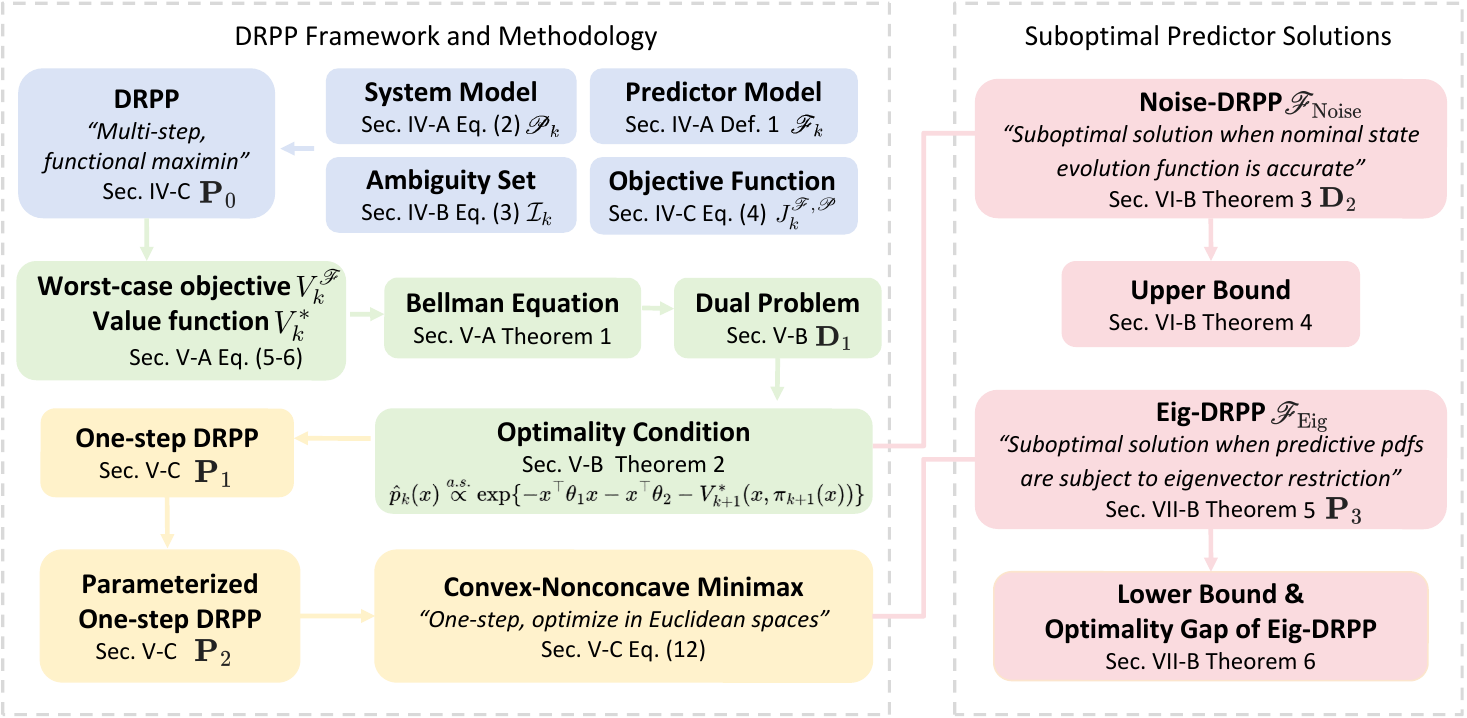}
    \caption{An illustration of the roadmap.}
    \label{fig:roadmap_2}
\end{figure}

The rest of this article is organized as follows. Section II provides a brief review of the related work and Section III introduces preliminaries on probabilistic prediction and DRO. In Section IV, we formulate the proposed DRPP problem. In Section V, we introduce the main methodology to address the main challenges in the DRPP. Then, we relax the constraints to solve two suboptimal DRPPs with upper and lower bounds in Section VI and Section VII. Next, we conduct elaborate experiments and provide an application in Section VIII, followed by a conclusion in Section IX. Please see Fig. \ref{fig:roadmap_2} for the overall structure of this paper.
% \begin{figure}[htb]
%     \centering
%     \includegraphics[width=0.9\linewidth]{pic/roadmap.pdf}
%     \caption{A roadmap of the paper's overall structure.}
%     \label{fig:roadmap}
% \end{figure}

\section{Related Work}
The study of DRPP for SDSs integrates ideas from fields including probabilistic prediction, control theory, minimax optimization, and DRO. In this section, we briefly summarize the most relevant work as follows.

\paragraph{Probabilistic prediction of SDS} According to how the predictive pdfs are approximated, existing literature in the control community can be classified into the following groups. 
First, approximation by the first several central moments. Particularly, there are many methods using the mean and covariance to describe the predictive distribution. Classical local approximation methods include the Taylor polynomial in the extended Kalman filter \cite{maybeckStochasticModelsEstimation1982, rothEfficientImplementationSecond2011}, the Stirling's interpolation \cite{norgaardNewDevelopmentsState2000} and the unscented transformation \cite{menegazSystematizationUnscentedKalman2015} as  representatives for derivative-free estimation methods \cite{simandlDerivativefreeEstimationMethods2009}.

Second, approximation by a weighted sum of basis functions. The Gaussian mixture model (GMM) \cite{sorensonRecursiveBayesianEstimation1971} approximates a pdf by a sum of weighted Gaussian distributions, and it can achieve any desired approximation accuracy if there are sufficiently large number of Gaussians \cite{alspachNonlinearBayesianEstimation1972}. Using GMM, the probabilistic prediction is equivalent to propagating the first two moments and weights of each basis \cite{chenMixtureKalmanFilters2000}.
The polynomial chaos expansion (PCE) \cite{wienerHomogeneousChaos1938} approximates a state evolution function by a sum of weighted polynomial basis functions, and these basis functions are orthogonal regarding the input pdf. In this way, the mean and covariance of the output can be conveniently approximated by the weights of PCE \cite{paulsonEfficientMethodStochastic2019}.

Third, approximation by a finite number of sampling points. The numerical solution using Monte Carlo (MC) samplings \cite{sarkkaBayesianFilteringSmoothing2023} can approximate any statistics of the output distribution if the sampling number is sufficiently high. For a further review of the modern MC methods in approximated probabilistic prediction, we recommend the review \cite{zhangModernMonteCarlo2021}.

Nevertheless, assuming the accuracy of a nominal model is restrictive, especially in the control engineering practice. As reviewed in \cite{landgrafProbabilisticPredictionMethods2023}, an important application of probabilistic prediction for SDSs is to formulate chance constraints \cite{farinaStochasticLinearModel2016} in stochastic model predictive control (SMPC) \cite{mesbahStochasticModelPredictive2016,mcallisterNonlinearStochasticModel2023}. In these situations, a misleading predictive pdf can lead to unsafe control policies. To take the distributional uncertainties into account, Coulson et al \cite{coulsonDistributionallyRobustChance2022} have proposed a novel distributionally robust data-enabled predictive control algorithm, Coppens and Patrinos \cite{coppensDataDrivenDistributionallyRobust2022a} develop a data-driven distributionally robust MPC scheme using generalized moment-based ambiguity sets. However, these works aim to robustify the control performance, while the problem of distributionally robust probabilistic prediction remains open.

\paragraph{Probabilistic prediction in machine learning} Taking a historical view of the development of PP in the statistics and machine learning communities, we can find that the research in the early phase also focuses on approximating the posterior pdf, especially from a Bayesian perspective. The advent of the scoring rule has provided a convenient scalar way to measure the performance of probabilistic predictions \cite{gneitingProbabilisticForecasting2014}. Then, the research focus is shifting toward training probabilistic predictors to optimize the empirical scores. Classical statistical methods include Bayesian statistical models \cite{robertBayesianChoice2007}, approximate Bayesian forecasting \cite{frazierApproximateBayesianForecasting2019}, quantile regression \cite{koenkerQuantileRegression402017}, etc. Recent machine learning algorithms include distributional regression \cite{schlosserDistributionalRegressionForests2019}, boosting methods such as NGBoost \cite{duanNGBoostNaturalGradient2020}, and autoregressive recurrent networks like DeepAR \cite{salinasDeepARProbabilisticForecasting2020}, etc. We recommend \cite{tyralisReviewPredictiveUncertainty2024} for a more detailed review of machine learning algorithms for probabilistic prediction.

\paragraph{Distributionally robust optimization}To guarantee the prediction performance for SDSs when the nominal model is not accurate, the primary task is to describe the ambiguity set for SDSs, which is the key concept in the DRO community. 
% DRO is a mathematical framework for decision-making under uncertainty. Moment-based ambiguity sets are among the most widely studied in DRO, leveraging partial statistical information (e.g., mean, covariance, or higher moments) to define a set of possible distributions.
In the seminal work \cite{scarfMinmaxSolutionInventory1957}, the Chebyshev ambiguity set is defined, where only the mean and covariance are known. Using this ambiguity set, several extensions are further studied in \cite{gallegoDistributionFreeNewsboy1993}. Next, Delage and Ye \cite{delageDistributionallyRobustOptimization2010} generalize the Chebyshev ambiguity set to the conic moment-based ambiguity set, which allows the mean and covariance matrix to be also unknown. Shapiro and Pichler \cite{shapiroConditionalDistributionallyRobust2024} further extend the conic moment-based ambiguity set to the conditional version, which is suitable for multistage decision-making problems.
Overall, viewing from the outer optimization, DRO still optimizes in an Euclidean space, while DRPP optimizes in a functional space of all pdfs. Since classical tools in DRO cannot be directly applied to solve DRPP, new methods are needed.

\section{Preliminaries}
In this paper, we denote random variables in \textbf{bold fonts} to distinguish them from constant variables. Given a random variable $\boldsymbol{x}$ taking values in $\mathcal{X}$, we denote its probability measure as $P_{\boldsymbol{x}}(\cdot)$ and its pdf as $p_{\boldsymbol{x}}(\cdot)$. Condsider another random variable $\boldsymbol{y}$, we deonte the conditional probability measure of $\boldsymbol{y}$ given $\boldsymbol{x}$ as $P_{\boldsymbol{x}\mid\boldsymbol{y}}(\cdot \mid \cdot)$ and the conditional pdf as $p_{\boldsymbol{x}\mid\boldsymbol{y}}(\cdot \mid \cdot)$. Let $\mathcal{P}(\mathcal{X})$ be the space of probability measures on $\mathcal{X}$, equipped with the standard weak*-topology. Let $\mathcal{P}_2(\mathcal{X})\subset\mathcal{P}(\mathcal{X})$ be the space of probability measures with finite second moments. Let $\mathcal{M}_{+}(\mathcal{X})$ be the set of positive measures on $\mathcal{X}$.
Given a measurable map $h:\mathcal{X}\to \mathbb{R}$, the expectation of $h(\boldsymbol{x})$ is denoted as $\mathbb{E}_{P_{\boldsymbol{x}}}h(\boldsymbol{x})$ or $\mathbb{E}_{x\sim p_{\boldsymbol{x}}} h(x)$. We also denote a sequence $\{(\cdot)_k\}_{k=1}^T$ by $(\cdot)_{1:T}$. We denote the set of symmetric matrices with dimension $n$ as $S^n$, whose subset of positive semi-definite matrices is $S_{+}^{n}$. For a function $f:\mathcal{X} \to \mathbb{R}$, the operator norm  $\|f\|_{2}:=\inf \{c \geq 0:\|f(v)\|_2 \leq c\|v\|_2 \;\forall v \in \mathbb{R}^n\}$. For $\rho \in \mathcal{P}(\mathcal{X})$, we denote $\langle \rho, f\rangle:= \int_{\mathcal{X}} f(x) \mathrm{d}\rho(x)$. Given $A, B\in\mathbb{R}^{n\times n}$, we use $A \succeq B$ to indicate $A-B \in S_{+}^{n}$, and we denote $\langle A, B \rangle := \operatorname{tr}(A^\top B)$.

\subsection{Probabilistic Prediction}
Given a random variable $\boldsymbol{x}$ taking value on $\mathcal{X}$, consider a predictive space, $\mathcal{F}$, of probability measures on $\mathcal{X}$ that admits a Lebesgue density, $\hat{p}_{\boldsymbol{x}}$, for each element. After the value of $\boldsymbol{x}$ is materialized as $x$, a scoring rule,
\begin{equation*}
    \mathcal{S}(\hat{p}_{\boldsymbol{x}}, x): \mathcal{F} \times \mathcal{X} \to \mathbb{R},
\end{equation*}
assigns a numerical score $\mathcal{S}(\hat{p}_{\boldsymbol{x}}, x)$ to measure the quality of the predictive distribution $\hat{p}_{\boldsymbol{x}}$ on the realized value $x$. Notice that the scoring rule $\mathcal{S}(\hat{p}_{\boldsymbol{x}}, \boldsymbol{x})$ is a random variable because it depends on the realization of $\boldsymbol{x}$.
% \textcolor{blue}{The expected scoring rule of $\mathcal{S}$ is denoted by $\mathbb{E}_{x\sim p_{\boldsymbol{x}}} \mathcal{S}(\hat{p}_{\boldsymbol{x}}, x)$ or $\mathbb{E}\mathcal{S}(\hat{p}_{\boldsymbol{x}}, \boldsymbol{x})$.}
A scoring rule $\mathcal{S}$ is proper with respect to the predictive space $\mathcal{F}$ if $\mathbb{E}\mathcal{S}(\hat{p}_{\boldsymbol{x}}, \boldsymbol{x}) \leq \mathbb{E}\mathcal{S}(p_{\boldsymbol{x}}, \boldsymbol{x})$ holds for all $\hat{p}_{\boldsymbol{x}}, p_{\boldsymbol{x}} \in \mathcal{F}$. It is strictly proper if the equality holds only when $\hat{p}_{\boldsymbol{x}}$ equals $ p_{\boldsymbol{x}}$ almost everywhere.
For example, the logarithm score, 
\[\mathcal{L}(\hat{p}_{\boldsymbol{x}}, x):= \log \hat{p}_{\boldsymbol{x}}(x),\] is one of the most celebrated scoring rules for being essentially the only local proper scoring rule up to equivalence \cite{parryProperLocalScoring2012a}.

\subsection{Distributionally Robust Optimization}
Traditional stochastic optimization approaches assume a known probability distribution $P$ and solve problems
\[
\inf_{x \in \mathcal{X}} \mathbb{E}_{\xi\sim P}[h(x, \xi)],
\]
where $x \in \mathcal{X}$ represents the decision variable,
$\xi$ is sampled from the distribution $P$, $h(x, \xi)$ is the objective function.
However, in many real-world problems, the true probability distribution of the uncertain parameters is unknown or partially known. 
DRO addresses this issue by considering a family of distributions, known as an ambiguity set, rather than a single distribution. The goal is to find a solution that performs well for the worst-case distribution within this set. A general DRO problem is formulated as:
\[
\inf_{x \in \mathcal{X}} \sup_{Q \in \mathcal{D}} \mathbb{E}_{\xi \sim Q}[h(x, \xi)],
\]
where $\mathcal{D}$ is the ambiguity set, a set of probability distributions that are close to the nominal distribution $P$ in some sense.

The ambiguity set $\mathcal{D}$ can be defined in various ways. One of the most popular ones is the moment-based sets, where distributions are constrained by their moments (e.g., mean, variance). For example, the celebrated conic moment-based ambiguity set is used in the seminal work \cite{delageDistributionallyRobustOptimization2010} such that
\begin{equation}\label{eq:classical_conic_moment-based_abg}
    \begin{aligned}
        &\mathcal{D}_1\left(\mathcal{X}, \mu_0, \Sigma_0, \gamma_{1}, \gamma_{2}\right)\\
        :=&\left\{\begin{array}{l|l}
    p_{\boldsymbol{\xi}} & \begin{array}{l}
    \mathbb{P}(\boldsymbol{\xi} \in \mathcal{X})=1 \\
    \left(\mathbb{E}[\boldsymbol{\xi}]-\mu_0\right)^{\top} \Sigma_0^{-1}\left(\mathbb{E}[\boldsymbol{\xi}]-\mu_0\right) \leq \gamma_{1} \\
    \mathbb{E}\left[\left(\boldsymbol{\xi}-\mu_0\right)\left(\boldsymbol{\xi}-\mu_0\right)^{\top}\right] \preceq \gamma_{2} \Sigma_0
    \end{array}
    \end{array}\right\},
    \end{aligned}
\end{equation}
where $\mathcal{X}$ is the nominal support, $\mu_0$ and $\Sigma_0$ are the nominal mean and covariance, $\gamma_{1}, \gamma_{2} \in \mathbb{R}_{+}$ quantifies the uncertainties on the mean and covariance.
A thorough review of popular ambiguity sets is provided in \cite[Sec. 3]{linDistributionallyRobustOptimization2022}.

\section{Problem Formulation}
To begin with, we define the system and probabilistic predictor model, respectively. Next, we generalize the classical conic moment-based ambiguity set \eqref{eq:classical_conic_moment-based_abg} for SDSs such that one's uncertainties of both state evolution functions and system noises are jointly described. Finally, we formulate the DRPP problem as a multistep functional maximin optimization.

\subsection{System and Predictor Model}
\textbf{System model ($\mathcal{X}, \mathcal{U}, \mathscr{P}, \pi, T$).} Consider a discrete-time SDS with horizon $T\in\mathbb{Z}_{+}$, denoted as $\mathscr{P}$, whose dynamics (or kernel) at time step $k\in\{0,\ldots, T-1\}$ is
\begin{equation}\label{eq:sds}
    \mathscr{P}_k: \boldsymbol{x}_{k+1} =f_k(\boldsymbol{x}_{k}, \boldsymbol{u}_k) + \boldsymbol{w}_{k},
\end{equation}
where $\mathscr{P}_k$ denotes the system kernel at time step $k$, $\boldsymbol{x}_k$ is the system state taking values in the state space $\mathcal{X} := \mathbb{R}^{d_x}$, the initial state is given as $x_0$.
$\boldsymbol{u}_k$ is the control taking values in the input space $\mathcal{U}:=\mathbb{R}^{d_u}$, which comes from a deterministic state-feedback policy $\pi_k: \mathcal{X} \to \mathcal{U}$, i.e., $\boldsymbol{u}_k = \pi_k(\boldsymbol{x}_k)$. The subsequent problem formulation thus evaluates predictive performance conditioned on a fixed policy $\pi$.
We denote the state-control pair $(\boldsymbol{x}_k,\boldsymbol{u}_k)$ as $\boldsymbol{z}_k \in \mathcal{Z} := \mathcal{X} \times \mathcal{U}$, then $f_k: \mathcal{Z} \to \mathcal{X}$ is the state evolution function.
$\boldsymbol{w}_k$ is an exogenous noise vector taking values in $\mathbb{R}^{d_x}$.

\begin{assumption}[SDS]
    At each time step $k\in\{0,\ldots, T-1\}$,
    \begin{itemize}
        \item $f_k$ has finite operator norms, i.e., $\|f_k\|_2$ is bounded, which is satisfied by many nonlinear systems in practice.
        \item the noise $\boldsymbol{w}_k$ has a bounded second moment, i.e., $P_{\boldsymbol{w}_k} \in \mathcal{P}_2(\mathbb{R}^{d_x})$, and $\boldsymbol{w}_{0:T-1}$ are independent but not necessarily identically distributed.
    \end{itemize}
\end{assumption}

\textbf{Predictor model ($\mathscr{F}$).} Starting from the system's initial state $x_0$, a probabilistic predictor keeps observing the states and the control inputs, aiming to predict the conditional pdfs of the next state.
\begin{definition}[Predictive space]\label{def:predictive_space}
    Let $\mathcal{F}\subset\mathcal{P}(\mathcal{X})$ be the set of continuous positive pdfs on $\mathcal{X}$ with tail envelope of order $2$, i.e.,
    \begin{equation*}
        \begin{aligned}
            \mathcal{F} :=& 
            \left\{\begin{array}{l|l}
                \!\!\!\hat{p} & 
                \!\begin{array}{l}
                \hat{p}\text{ is a strictly positive pdf on }\mathcal{X},\\
                \exists\; C>0\text{ s.t.}|\log \hat p(x)| \le C(1+\|x\|_2^2)\; \forall x\in\mathcal{X}
                \end{array}
            \end{array}\!\!\!\!\!\right\}.
        \end{aligned} 
    \end{equation*}
\end{definition}
A probabilistic predictor policy for $\mathscr{P}$ is a sequence $\mathscr{F} = (\mathscr{F}_0, \ldots, \mathscr{F}_{T-1})$ where each $\mathscr{F}_k: \mathcal{Z} \to \mathcal{F}$ is a measurable map from the previous state and control input to a predictive conditional pdf for the next state. Let $\mathfrak{F}$ be the collection of such predictor policies.

\begin{assumption}[Predictor's prior information]\label{assumption:prior-info}
    At each time step $k\in\{0,\ldots, T-1\}$, the probabilistic predictor has only access to the following information:
    \begin{itemize}
        \item the state trajectory $x_{0:k}$, control input sequence $u_{0:k}$,
        \item a nominal state evolution function $\bar{f}_k:\mathcal{Z} \to \mathcal{X}$,
        \item a nominal noise mean $\bar{\mu}_k\in\mathcal{U}$, and a nominal noise covariance $\bar{\Sigma}_k\in S_{+}^{d_x}$.
    \end{itemize}
\end{assumption}
The predictor's prior information usually does not align with the ground truth. Therefore, to further quantify the uncertainty between the nominal model and the real model, we need to define an ambiguity set of SDSs.

\subsection{Ambiguity Set}
The crucial spirit of an ambiguity set is that although a nominal model rarely identifies with the ground truth, one can still have additional information or belief that the real model is located around the nominal model.
For example, an upper bound for the operator norm $\|f_k-\bar{f}_k\|_2$ may be available when the nominal function is derived from a system identification procedure. 
Additionally, estimations and confidence regions (CRs) of the mean and covariance of system noise may be statistically inferred from historical data. Even if no supportive evidence exists to construct a reasonable ambiguity set, one can still subjectively select an ambiguity set.

Jointly considering the uncertainties of both state evolution functions and noises, we define the conditional conic moment-based ambiguity sets for $\mathscr{P}$.
\begin{definition}[Conditional conic moment-based ambiguity set]\label{def:ambiguity_set}
    For each time step $k\in\{0,\ldots,T-1\}$ and state-control pair $z$, the conditional conic moment-based ambiguity set is a subset of $\mathcal{P}_2(\mathcal{X})$ defined as
    \begin{equation}\label{eq:ambiguity-set}
        \begin{aligned}
            &\quad\mathcal{I}_k\left(z \mid \bar{f}_k, \bar{\mu}_k, \bar{\Sigma}_k, \gamma_{0}, \gamma_{1}, \gamma_{2}, \gamma_{3}\right)=:\\
            &\left\{\!\!\!\begin{array}{l|l}
            P_{\boldsymbol{x}_{k\!+\!1}\mid\boldsymbol{z}_k}(\cdot\!\mid\! z)\! &\!\!\! \begin{array}{l}
            \boldsymbol{x}_{k+1}=f_k(z)+\boldsymbol{w}_k,\\
            \|f_k(z)-\bar{f}_k(z)\|_{2}^2 \leq \gamma_{0}(z)\\
            \|\mathbb{E}(\boldsymbol{w}_k)-\bar{\mu}_k\|_{\bar{\Sigma}_k^{-1}}^2 \leq \gamma_{1} \\
            \gamma_{3}\bar{\Sigma}_k \!\preceq\! \mathbb{E}\left(\!\boldsymbol{w}_k\!-\!\bar{\mu}_k\right)\!\left(\boldsymbol{w}_k\!-\!\bar{\mu}_k\right)^{\!\top} \!\!\preceq\! \gamma_{2} \bar{\Sigma}_k
            \end{array}
            \end{array}\!\!\!\!\!\right\},
        \end{aligned}
    \end{equation}
    where $\bar{f}_k, \bar{\mu}_k, \bar{\Sigma}_k$ are defined in Assumption \ref{assumption:prior-info}, and predictor's uncertainties are quantified by $\gamma_{0}: \mathcal{Z} \to \mathbb{R}_{+}$ and $\gamma_j \in \mathbb{R}_{+}$ for $j=1,2,3$. 
    When there is no confusion about the parameters that define an ambiguity set, we use $\mathcal{I}_k(z)$ for simplification. 
\end{definition}
The conic moment-based ambiguity set covers a wide family of distributions with bounded second moments. This choice not only facilitates tractable dual reformulations but also provides a statistically interpretable description of moment uncertainty, leading to guaranteed robustness without assuming a specific parametric form.
\begin{assumption}[Regularity]\label{ass:regularity}
    At each time step $k\in\{0,\ldots,T-1\}$, conditioned on the state-control pair $\boldsymbol{z}_k=z_k$, the system chooses a conditional probability measure $P_k(\cdot \mid z_k) \in \mathcal{I}_k(z_k)$ indepentdent of the previous choices.
\end{assumption}
Let $\mathfrak{P}$ be a collection of SDSs subject to the conditional conic moment-based ambiguity set and regularity assumption.

\subsection{Problem in Interest}
\textbf{Objective function ($J_k^{\mathscr{F}, \mathscr{P}}$)}. Starting from the time step $k\in\{0,\ldots,T-1\}$, given an initial state-control pair $z\in\mathcal{Z}$, a predictor policy $\mathscr{F}\in\mathfrak{F}$, and a system kernel $\mathscr{P}\in\mathfrak{P}$, the prediction performance of $\mathscr{F}$ is characterized by the expectation of cumulative log-score over the state trajectory. We define it by the objective function as follows:
\begin{equation}\label{eq:obj}
    J_k^{\mathscr{F}, \mathscr{P}}(z) := \mathbb{E}_{\mathscr{P}_{k:T-1}}\!\left[ \sum_{t=k}^{T-1} \mathcal{L}\left( \mathscr{F}_t(\boldsymbol{z}_t), \boldsymbol{x}_{t+1} \right) \bigg| \boldsymbol{z}_k = z\right],  
\end{equation}
where the notation $\mathbb{E}_{\mathscr{P}}$ indicates that the expectation is taken over the random trajectory ($\boldsymbol{x}_k,\dots,\boldsymbol{x}_T$) generated by the kernel sequence $\mathscr{P}_{k:T-1}$.

\textbf{Distributionally robust probabilistic prediction ($\mathbf{P}_0$).} Given a set of probabilistic predictors $\mathfrak{F}$ and a set of system kernels $\mathfrak{P}$ for the SDS $\mathscr{P}$ with an initial state-control pair $z_0$, we are interested in finding probabilistic predictors in $\mathfrak{F}$ that are distributionally robust with respect to $\mathfrak{P}$. In other words, a distributionally robust probabilistic predictor should have a worst-case performance guarantee for any system kernel in $\mathfrak{P}$. To this end, a multistep functional maximin optimization is formulated as follows:
\begin{equation*}
    (\mathbf{P}_0): \sup_{\mathscr{F}\in\mathfrak{F}}\; \inf_{\mathscr{P}\in\mathcal{\mathfrak{P}}} \; J_0^{\mathscr{F}, \mathscr{P}}(z_0).
\end{equation*}

\section{Main Methodology}
To begin with, we derive the Bellman equation of DRPP based on the principle of dynamic programming, which shows that solving the DRPP problem $\mathbf{P}_0$ is equivalent to solving a sequence of Bellman equations in a backward manner.
Then, we canonicalize the ambiguity set of SDS and transform the inner optimization into an infinite-dimensional conic linear program. 
Next, we take the dual form of inner minimization with a strong duality guarantee. A necessary optimality condition is developed to handle the computational intractability of functional optimization. 
Exploiting this optimality condition, we equivalently transform the one-step DRPP in function spaces into a convex-nonconcave minimax problem in Euclidean spaces.
Finally, a discussion on the computational complexity is provided.

\subsection{Bellman Equation}
The concept of robust value function originates from the study of robust Markov decision processes (RMDPs), where the optimal worst-case performance of a control problem is quantified. Here, we apply the same idea to $\mathbf{P}_0$. 
Starting from the time step $k\in\{0,\ldots,T-1\}$, given an initial state-control pair $z\in\mathcal{Z}$ and a predictor policy $\mathscr{F}\in\mathfrak{F}$, the worst-case objective function for $\mathscr{F}$ is defined as 
\begin{equation}\label{eq:worst-obj}
    V^\mathscr{F}_k(z) := \inf_{\mathscr{P} \in \mathfrak{P}} J_k^{\mathscr{F}, \mathscr{P}}(z). 
\end{equation}
Then, the robust value function is defined as the optimal worst-case objective function, i.e.,
\begin{equation}\label{eq:robust-value}
   V_k^*(z) := \sup_{\mathscr{F}\in\mathfrak{F}} V_k^{\mathscr{F}}(z).
\end{equation}

Leveraging the principle of dynamic programming, one can derive the following recursive equations for $V_k^*$.
\begin{theorem}[Bellman equation]\label{thm:DPP}
    The set of robust value functions $\{V_k^*, k = 0, \ldots, T\}$ satisfies the following Bellman equation:
    for any state-control pair $z\in\mathcal{Z}$, $V^*_T(z) = 0$, and for $k\in\{0,\ldots,T-1\}$, there is
    \begin{equation}\label{eq:bellman}
        V^*_k(z) \!=\! \sup_{\hat{p}_k \in \mathcal{F}} \inf_{\rho_k \in \mathcal{I}_k(z)} \!\int_{\mathcal{X}} \!\mathcal{L}(\hat{p}_k, x) + V^*_{k+1}(x, \pi_{k+1}(x)) \mathrm{d} \rho_k(x).
    \end{equation}
\end{theorem}
\begin{proof}
    Please see Appendix \ref{app:thm:DPP}.
\end{proof}
\begin{remark}
    The Bellman equation implies that different control policies can lead to different robust value functions even when the system kernel family $\mathfrak{P}$ and predictor family $\mathfrak{F}$ are fixed. 
    Thus, the DRPP $\mathbf{P}_0$ can be extended to jointly optimize the control policy and the predictor policy in future work.
\end{remark}
Theorem \ref{thm:DPP} shows that solving the DRPP problem $\mathbf{P}_0$ is equivalent to solving the Bellman equation \eqref{eq:bellman}. This result generalizes the classical Bellman equation for RMDPs \cite{iyengarRobustDynamicProgramming2005} to the DRPP setting.
 
Let $\nu_k = f_k(z)$ and $\bar\nu_k = \bar f_k(z)$, the ambiguity constraint of $f_k$ is equivalent to describing the distance between $\nu_k$ and $\bar\nu_k$, i.e., $\|\nu_k-\bar\nu_k\|_2^2=\|f_k(z)- \bar f_k(z)\|_2^2\leq \gamma_{0}(z)$.
To deal with the remaining two ambiguity constraints on the first two moments of $\boldsymbol{w}_k$, let $\mu_k = \mathbb{E}\boldsymbol{w}_k$, and the r.h.s. of \eqref{eq:bellman} is transformed as
\begin{equation}\label{eq:p_0_cannolized}
    \begin{aligned}
        &\sup_{\hat{p}_k\in\mathcal{F}}\inf_{\nu_k, \mu_k, \rho_k} \int_{\mathcal{X}}\!\! \mathcal{L}(\hat{p}_k, x) + V^*_{k+1}\left(x, \pi_{k+1}(x)\right) \mathrm{d} \rho_k(x)\\
        &\text{s.t.} \left\{
        \begin{aligned}
            & \rho_k\in\mathcal{M}_{+}(\mathcal{X}),\mathbb{E}_{\rho_k}[1]=1, \mathbb{E}_{x\sim \rho_k}[x]=\mu_k+\nu_k  \\
            & \gamma_{3}\bar{\Sigma}_k \!\preceq\! \mathbb{E}_{x\sim \rho_k}\!\left[\left(x\!-\!\nu_k\!-\!\bar{\mu}_k\right)\left(x\!-\!\nu_k\!-\!\bar{\mu}_k\right)^{\top}\right] \!\preceq\! \gamma_{2} \bar{\Sigma}_k \\
            & \left[\begin{array}{cc}
            \bar{\Sigma}_k & \left(\mu_k-\bar{\mu}_k\right) \\
            \left(\mu_k-\bar{\mu}_k\right)^{\top} & \gamma_{1}
            \end{array}\right] \succeq 0\\ &\left[\begin{array}{cc}
            I & \left(\nu_k-\bar\nu_k\right) \\
            \left(\nu_k-\bar\nu_k\right)^{\top} & \gamma_{0}(z)
            \end{array}\right] \succeq 0,
        \end{aligned}
        \right.
    \end{aligned}
\end{equation}
where the inner problem is canonicalized into an infinite-dimensional conic linear program.

\subsection{Optimality Condition}
Directly solving the worst-case probability measure in the positive measure space $\mathcal{M}_{+}(\mathcal{X})$ is computationally intractable. As a first step towards dealing with this challenge, we leverage the dual analysis for \eqref{eq:p_0_cannolized}. If a strong duality holds, one can optimize the Lagrange multipliers in the dual form. For the inner optimization of \eqref{eq:p_0_cannolized}, a classical conclusion is that $\bar{\Sigma}_k \succ 0$ is a sufficient condition for strong duality to hold \cite[p. 5]{delageDistributionallyRobustOptimization2010}. Let $r\in\mathbb{R},q\in\mathbb{R}^{d_x}, Q_i\in S_{+}^{d_x}, \kappa_i = \{P_i \in S_{+}^{d_x}, p_i\in \mathbb{R}^{d_x}, s_i\in\mathbb{R}\}$ for $i=1,2$ be the Lagrange multipliers, the dual form of \eqref{eq:p_0_cannolized} is
\begin{equation*}\label{eq:p_0_cannolized-dual-join}
    \begin{aligned}        (\mathbf{D}_1):&\sup_{\hat{p}_k\in\mathcal{F},r,q,Q_1,Q_2,\kappa_1,\kappa_2}  G(r,q,Q_1,Q_2,\kappa_1,\kappa_2)\\          
        &\text{s.t.} \left\{
        \begin{aligned}
            &x^\top(Q_1-Q_2)x + x^\top q + r + \log \hat{p}_k(x)\\
            &+ V_{k+1}^*(x,\pi_{k+1}(x))\geq 0 \; \forall x \in \mathcal{X}\\
            % & q+2p_1+2(Q_1-Q_2)(\nu_k+\bar{\mu}_k) = 0\\
            &\left[\begin{array}{cc}
                P_1 & p_1\\
                p_1^\top & s_1
            \end{array}\right]\succeq 0, \left[\begin{array}{cc}
                P_2 & p_2\\
                p_2^\top & s_2
            \end{array}\right]\succeq 0.
        \end{aligned}
        \right.
    \end{aligned}
\end{equation*}
where $G(r,q,Q_1,Q_2,\kappa_1,\kappa_2) \!=\! - r \!+\! \bar{\mu}_k^\top\!(Q_1\!-\!Q_2)\bar{\mu}_k- (\gamma_{2}Q_1 - \gamma_{3}Q_2+P_1)\cdot\bar{\Sigma}_k- P_2\cdot I + 2p_1^\top\bar{\mu}_k- s_1\gamma_1+2p_2^\top\bar\nu_k - s_2\gamma_{0}(z) + g(q,p_1,p_2,Q_1,Q_2)$, and $g(q,p_1,p_2,Q_1,Q_2) = \inf_{\nu_k} - (q+2p_2)^\top\nu_k - \nu_k^\top(Q_1-Q_2)\nu_k\text{ s.t. } q+2p_1+2(Q_1-Q_2)(\nu_k+\bar{\mu}_k) = 0.$ Please see Appendix \ref{app:lem:lower-dual-problem} for a detailed derivation.

Unfortunately, we have to deal with an infinite number of constraints indexed by $x\in\mathcal{X}$ in problem $\mathbf{D}_1$. Although the dual problem cannot lead to a direct solution, it helps to develop an optimality condition for the probabilistic predictors in the following theorem.
\begin{theorem}[Optimality condition]\label{thm:opt-condition}
    A necessary optimality condition for the DRPP $\mathbf{P}_0$ is that $\hat{p}_k$ is belongs to the following exponential family almost everywhere on $\mathbb{R}^{d_x}$, 
    \begin{equation}\label{eq:exp-family}
        \hat{p}_k(x) \overset{a.s.}{\propto} \exp\{-x^\top \theta_1 x - x^\top \theta_2 - V_{k+1}^*(x, \pi_{k+1}(x))\},
    \end{equation}
    where $\propto$ means the right-hand side gives the density up to a normalizing constant, $\theta_1 \in S^{d_x}$, $\theta_2 \in \mathbb{R}^{d_x}$, and $k\in\{0,\ldots,T-1\}$. Particularly for $k=T-1$, the optimial $\hat{p}_{T-1}$ is subject to a Gaussian distribution almost everywhere, i.e., $\exists\;\hat{\mu}_k\in\mathbb{R}^{d_x}$ and $\hat{\Sigma}_k \in S^{d_x}_{+}$ such that
    \begin{equation}\label{eq:gaussian}
        \hat{p}_k \overset{\text{a.s.}}{\sim} \mathcal{N}\left(\hat{\mu}_k, \hat{\Sigma}_k\right).
    \end{equation}
\end{theorem}
\begin{proof}
    Please see Appendix \ref{app:thm:opt-condition}.
\end{proof}
% \textcolor{blue}{\begin{remark}
%     The conclusion that the optimal predictor allows for a Gaussian parameterization is derived from the structural characteristics of conic moment-based constraints and the logarithm scoring rule. Intuitively, moment constraints only characterize a distribution by means of its mean and covariance; under the objective of the log-score, these moment constraints compel the worst-case optimal member of the admissible family to take a Gaussian form.
% \end{remark}}

Theorem \ref{thm:opt-condition} has reduced the optimization complexity from the original function space $\mathcal{F}$ to a specific exponential distribution family. For $k=T-1$, the space is further reduced to $S_{+}^{d_x} \times \mathbb{R}^{d_x}$, 
where $\mathbf{D}_1$ can be transformed into a finite-dimensional optimization problem as follows
\begin{equation}\label{eq:p_0_cannolized-dual-join-parametrize}
    \begin{aligned}
        &\max_{\hat{\mu}_k,\hat{\Sigma}_k,Q_1,Q_2,\kappa_1,\kappa_2} -\frac{1}{2}d_x\log(2\pi) \!+\! \frac{1}{2}\log\operatorname{det}(\hat{\Sigma}_k^{-1})\\
        &- (\gamma_{2}Q_1\!-\!\gamma_{3}Q_2)\!\cdot\!\bar{\Sigma}_k-P_1\!\cdot\!\bar{\Sigma}_k - P_2\!\cdot\! I - s_1\gamma_1\\
        &- s_2\gamma_{0}(z_)^2 + 2p_1^\top\!\hat{\Sigma}_k (2p_2\!-\!p_1)- 2p_2^\top(\hat{\mu}_k\!-\!\bar{\mu}_k\!-\!\bar{\nu}_k)\\      
        &\text{s.t.} \left\{
        \begin{aligned}
            &Q_1 \succeq 0, Q_2 \succeq 0, Q_1 - Q_2 = \frac{1}{2}\hat{\Sigma}_k^{-1}, \hat{\Sigma}_k \succ 0\\
            &\left[\begin{array}{cc}
                        P_1 & p_1\\
                        p_1^\top & s_1
                    \end{array}\right]\succeq 0, \left[\begin{array}{cc}
                        P_2 & p_2\\
                        p_2^\top & s_2
                    \end{array}\right]\succeq 0.
        \end{aligned}
        \right.
    \end{aligned}
\end{equation}
Please see Appendix \ref{app:use_optimality_condition} for a detailed derivation.
However, this problem is now a nonconvex semidefinite optimization with respect to $\hat{\Sigma}_k$, solving a globally optimal solution is NP-hard in general.
Even worse, the lack of a closed-form solution of $V^*_{T-1}$ makes it impossible to recursively parameterize and solve $\hat{p}_{k}$ backwardly based on the Bellman equation \eqref{eq:bellman}.

\subsection{One-step DRPP}
At time step $k\in\{0,\ldots, T-1\}$, the Bellman equation \eqref{eq:bellman} implies that the optimal prediction policy  maximizes the one-step score $\mathcal{L}(\hat{p}_k, x)$ plus the expected best future cumulative scores $V_{k+1}^*(x, \pi_{k+1}(x))$ from the next state $x$ onward. If one does not take into account the future scores, Theorem \ref{thm:opt-condition} shows that the optimal predictive pdf can be parametrized by a Gaussian distribution. Given $z\in\mathcal{Z}$, consider the following \textbf{one-step DRPP} problem:
\begin{equation*}
    (\mathbf{P}_1): \sup_{\hat{p}_k\in\mathcal{F}}\inf_{\rho_k\in\mathcal{I}_k(z)}\int_{\mathcal{X}} \mathcal{L}(\hat{p}_k, x) \mathrm{d}\rho_k(x),
\end{equation*} 
we substitute $\hat{p}_k$ by \eqref{eq:gaussian} and the objective follows as
\begin{equation*}
    \begin{aligned}
        &\mathbb{E}_{x\sim \rho_k} \!\!-\!\frac{1}{2}[d_x\log(2\pi)
        \!+\!\log\operatorname{det}(\hat{\Sigma}_k) \!+\! (x\!-\!\hat{\mu}_k)\hat{\Sigma}_k^{-1}(x\!-\!\hat{\mu}_k)\!^\top].
    \end{aligned}
\end{equation*}
Next, notice that $\mathbb{E}_{x\sim \rho_k}[x]=\mu_k+\nu_k$ and $\mathbb{E}_{x\sim \rho_k}(x-\nu_k\!-\mu_k)(x-\nu_k-\mu_k)^\top=\Sigma_k$, one has $\mathbb{E}_{x\sim \rho_k} (x-\hat{\mu}_k)\hat{\Sigma}_k^{-1}(x-\hat{\mu}_k)^\top$ equals to $[\Sigma_k + (\mu_k+\nu_k-\hat{\mu}_k)(\mu_k+\nu_k-\hat{\mu}_k)^\top]\cdot \hat{\Sigma}_k^{-1}$. The functional minimax $\mathbf{P}_1$ is then equivalent to a minimax optimization in Euclidean spaces as follows
\begin{equation*}\label{eq:one-step-DRPP-general-parameterized-1}
    \begin{aligned}
        (\mathbf{P}_2):&\min_{\hat{\Sigma}_k,\hat{\mu}_k}\max_{\Sigma_k,\mu_k,\nu_k} -\log\operatorname{det}(\hat{\Sigma}_k^{-1})\\
        &\qquad+\left[\Sigma_k\!+\!(\mu_k\!+\!\nu_k\!-\!\hat{\mu}_k)(\mu_k\!+\!\nu_k\!-\!\hat{\mu}_k)^\top\right]\!\cdot\! \hat{\Sigma}_k^{-1} \\
        & \text{s.t.}\left\{
        \begin{aligned}
            & \|\nu_k-\bar{\nu}_k\|_2^2 \leq \gamma_{0}(z),\|\mu_k-\bar{\mu}_k\|_{\bar{\Sigma}_k^{-1}}^2 \leq \gamma_1\\
            & \gamma_{3}\bar{\Sigma}_k\preceq\Sigma_k + (\mu_k-\bar{\mu}_k)(\mu_k-\bar{\mu}_k)^\top \preceq \gamma_{2}\bar{\Sigma}_k.
        \end{aligned}
        \right.
    \end{aligned}
\end{equation*}

\begin{lemma}\label{lem:primal-solution-1}
    Given any $\mu_k$, the one-step DRPP $\mathbf{P}_2$ has explicit solutions for $\Sigma_k$ and $\hat{\mu}_k$ as\\
    i) $\Sigma_k^{*} = \gamma_{2}\bar{\Sigma}_k - (\mu_k-\bar{\mu}_k)(\mu_k-\bar{\mu}_k)^\top$.\\
    ii) $\hat{\mu}^{*}_k = \bar{\mu}_k + \bar{\nu}_k$.
\end{lemma}
\begin{proof}
    Please see Appendix \ref{app:lem:primal-solution-1}.
\end{proof}

Substituting the expressions of $\Sigma_k^{*}$ and $\hat{\mu}_k^{*}$ into $\mathbf{P}_2$, it is equivalent to a convex-nonconcave minimax optimization:
\begin{equation}\label{eq:qcqp}
    \begin{aligned}
    	&\min_{\hat{\Sigma}_k\succ 0} \max_{a,b}\; -\log\operatorname{det}(\hat{\Sigma}_k^{-1}) + \gamma_{2}\bar{\Sigma}_k \cdot \hat{\Sigma}_k^{-1} \\
        &\qquad\qquad+ \|a+b\|_{\hat{\Sigma}_k^{-1}}^2 - \|a\|_{\hat{\Sigma}_k^{-1}}^2 \\
    	&\text{s.t.}\;  \|a\|_{\bar{\Sigma}_k^{-1}}^2\leq \gamma_1, \|b\|_2^2 \leq \gamma_{0}(z),
    \end{aligned}
\end{equation}
where $a = \mu_k-\bar{\mu}_k, b = \nu_k-\bar{\nu}_k$. Notice that the inner maximization is an indefinite quadratic constrained quadratic programming (QCQP) problem, and there are still short of efficient algorithms to solve this kind of minimax optimization.

Given any $\hat{\Sigma}_k \succ 0$ and $b\neq 0$, one can first solve $a$ by considering the following quadratic constrained linear programming (QCLP) problem:
\begin{equation}\label{eq:QCLP}
    \max_{a} \|a+b\|_{\hat{\Sigma}_k^{-1}}^2 - \|a\|_{\hat{\Sigma}_k^{-1}}^2 \text{ s.t. }  \|a\|_{\bar{\Sigma}_k^{-1}}^2\leq \gamma_1. 
\end{equation}
\begin{lemma}\label{lem:QCLP}
    The solution of the QCLP problem \eqref{eq:QCLP} is
    \begin{equation*}
        a^{*} = \frac{\sqrt{\gamma_1}\bar\Sigma_k\hat{\Sigma}_k^{-1}b}{\|\bar\Sigma_k\hat{\Sigma}_k^{-1}b\|_{\bar{\Sigma}_k^{-1}}}.
    \end{equation*}
\end{lemma}
\begin{proof}
    % Please see Appendix \ref{app:lem:QCLP}.
    Using the KKT condition, there are
    $$\left\{\begin{aligned}
    &\partial_{a} \left[\|a+b\|_{\hat{\Sigma}_k^{-1}}^2 - \|a\|_{\hat{\Sigma}_k^{-1}}^2 - s_1(\|a\|_{\bar{\Sigma}_k^{-1}}^2 - \gamma_1)\right] = 0\\
    & s_1(\|a\|_{\bar{\Sigma}_k^{-1}}^2 - \gamma_1) = 0, s_1\geq 0\\
    \end{aligned}\right.$$
    The first condition reveals that $s_1 \neq 0$ and $a = s_1^{-1}\bar{\Sigma}_k\hat{\Sigma}_k^{-1}b$, thus the second condition leads to $\|a\|_{\bar{\Sigma}_k^{-1}}^2 = \gamma_1$. Then, 
    \[s_1 = \frac{\|\bar\Sigma_k\hat{\Sigma}_k^{-1}b\|_{\bar{\Sigma}_k^{-1}}}{\|a\|_{\bar{\Sigma}_k^{-1}}}=\frac{\|\bar\Sigma_k\hat{\Sigma}_k^{-1}b\|_{\bar{\Sigma}_k^{-1}}}{\sqrt{\gamma_1}},\]
    and the proof is completed.
\end{proof}
Substituting the optimal $a^{*}$ into the problem, the inner optimization of \eqref{eq:qcqp} can be further reformulated as a convex maximization problem:
\begin{equation}\label{eq:convex-max}
    \max_{b} \; \alpha\|b\|_{A} + \|b\|_{B}^2 \text{ s.t. } \|b\|_2^2 \leq \gamma_{0}(z),
\end{equation}
where $A = \hat{\Sigma}_k^{-1}\bar{\Sigma}_k\hat{\Sigma}_k^{-1}$, $B=\hat{\Sigma}_k^{-1}$, and $\alpha = 2\sqrt{\gamma_1}$. 

\subsection{Complexity Analysis}
As we have shown, solving the one-step DRPP $\mathbf{P}_1$ is equivalent to finding a global minimax point (Stackelberg equilibrium) for a convex-nonconcave minimax problem where the outer minimization is a semidefinite programming and the inner problem is to maximize a non-smooth convex function subject to a quadratic constraint. It is known that convex maximization is NP-hard in very simple cases, such as quadratic maximization over a hypercube, and even verifying local optimality is NP-hard \cite{pardalosCheckingLocalOptimality1988a}. Thus, finding a global minimax point of a one-step DRPP problem is also NP-hard.

Even worse, finding a local minimax point of a general constrained minimax optimization is of fundamental hardness. Not only may the local minimax point not exist \cite{jinWhatLocalOptimality2020}, but any gradient-based algorithm needs exponentially many queries in the dimension and $\epsilon^{-1}$ to compute an $\epsilon$-approximate first-order local stationary point \cite{daskalakisComplexityConstrainedMinmax2021}. Finally, the local stationary point is not always guaranteed to be a local minimax point \cite{grimmerLimitingBehaviorsNonconvexNonconcave2022, fiezImplicitLearningDynamics2020a}.

In summary, globally finding a minimax point of a one-step DRPP is computationally intractable, and gradient-based algorithms aiming for local minimax points are limited by both convergence speed and performance guarantee. 
Rather than pursuing the optimal solutions of $\mathbf{P}_0$, a more reasonable goal is to derive suboptimal solutions for the one-step DRPP $\mathbf{P}_1$. Based on these suboptimal solutions, we can then derive upper and lower bounds for the robust value functions $V_k^*$.

\section{Value Function Upper Bound}
To get an upper bound of $V_k^{*}(z)$, one can either enlarge the feasible region of the outer maximization or shrink the feasible region of the inner minimization. In this paper, we choose to shrink the ambiguity set $\mathcal{I}_k(z)$.
If the first constraint is strengthened as $\nu_k = \bar{\nu}_k$, one is constrained to a smaller ambiguity subset of $\mathcal{I}_k(z)$. If a global minimax point can be solved for the relaxed one-step DRPP, one will get a suboptimal predictor and a performance upper bound. 
\subsection{Problem Reformulation}
If we shrink the ambiguity set of one-step DRPP $\mathbf{P}_1$ by assuming $\bar{\nu}_k=\nu_k$, there is $\boldsymbol{w}_k = \boldsymbol{x}_{k+1} - \bar{\nu}_k,$ which means that $\boldsymbol{w}_k$ can be precisely known after $\boldsymbol{x}_{k+1}$ been observed. 
Now that there is no uncertainty of the state evolution $\nu_k$, predicting the state's conditional probability measure $P_{\boldsymbol{x}_{k+1}\mid\boldsymbol{z}_k}$ is equivalent to predicting the noise' probability measure $P_{\boldsymbol{w}_{k}}$, and the one-step DRPP $\mathbf{P}_1$ is then reformulated as
\begin{equation}\label{eq:noise-DRPP}
    \begin{aligned}
        &\sup_{\hat{p}_k\in\mathcal{F}}\inf_{\mu_k, \tilde{\rho}_k} \; \int_{\mathcal{X}} \log [\hat{p}_k(w + \bar{\nu}_k)] \mathrm{d} \tilde{\rho}_k(w)\\
        &\text{s.t.} \left\{
        \begin{aligned}
            & \tilde{\rho}_k\in\mathcal{M}_{+}(\mathcal{X}),\mathbb{E}_{w\sim \tilde{\rho}_k}[1]=1, \mathbb{E}_{w\sim \tilde{\rho}_k}[w]=\mu_k  \\
            & \gamma_{3}\bar{\Sigma}_k \preceq \mathbb{E}_{w\sim \tilde{\rho}_k}\!\left[\left(w-\bar{\mu}_k\right)\left(w-\bar{\mu}_k\right)^\top\right] \preceq \gamma_{2} \bar{\Sigma}_k \\
            & \left[\begin{array}{cc}
            \bar{\Sigma}_k & \left(\mu_k-\bar{\mu}_k\right) \\
            \left(\mu_k-\bar{\mu}_k\right)^{\top} & \gamma_1
            \end{array}\right] \succeq 0.\\
        \end{aligned}
        \right.
    \end{aligned}
\end{equation}
With the outer maximization imposed on $\hat{p}_{k}$, we can get a reformulated one-step DRPP whose optimal worst-case value function is an upper bound of the objective of $\mathbf{P}_1$.

\subsection{Noise-DRPP}
Taking the dual form of the inner problem of \eqref{eq:noise-DRPP} and joining it with the outer maximization, we transform the primal maximin problem into the following maximization problem:
\begin{equation*}\label{eq:noise-DRPP-dual-join}
    \begin{aligned}
         (\mathbf{D}_2):&\sup_{\hat{p}_k\in\mathcal{F},r,q,Q_1,Q_2,P,p,s} G(r,q,Q_1,Q_2,P,p,s)\\          
        &\text{s.t.} \left\{
        \begin{aligned}
            &w^\top\!(Q_1\!-\!Q_2)w \!+\! w^\top q \!+\! r \!+\! \log \hat p_k(w\!+\!\bar\nu_k) \!\geq\! 0\\
            &\qquad\qquad\qquad\qquad\qquad\;\;\forall w \!\in\! \mathbb{R}^{d_x}\\
            &q + 2(Q_1-Q_2)\bar{\mu}_k + 2p = 0\\
            &Q_1 \succeq 0, Q_2 \succeq 0\\
            &\left[\begin{array}{cc}
                P & p\\
                p^\top & s
            \end{array}\right]\succeq 0,
        \end{aligned}
        \right.
    \end{aligned}
\end{equation*}
where $G(r,q,Q_1,Q_2,P,p,s) = - r + \bar{\mu}_k^\top(Q_1-Q_2)\bar{\mu}_k- (\gamma_{2}Q_1 - \gamma_{3}Q_2+P)\cdot\bar{\Sigma}_k+ 2p^\top\bar{\mu}_k- s\gamma_1$. Please see Appendix \ref{app:lem:upper-dual-problem} for a detailed derivation of the dual problem.

\begin{theorem}[Noise-DRPP]\label{thm:one-step-DRPP-mp}
    Given $\boldsymbol{z}_k = z \in \mathcal{Z}$ at time step $k\in\{0,\ldots,T-1\}$, if the probabilistic predictor has no ambiguity of the one-step state evolution, i.e., $\bar{\nu}_k = \nu_k$, the solution to the one-step DRPP $\mathbf{P}_1$ is:\\
    i) The optimal predictive pdf is 
    \begin{equation*}
        \hat{p}^{*}_k \sim \mathcal{N}\left(\bar{\nu}_k + \bar{\mu}_k, \;\gamma_{2}\bar{\Sigma}_k\right).
    \end{equation*}
    ii) The worst-case conditional measure $\rho_k^*$ belongs to the set
    \begin{equation*}
        \left\{\begin{array}{c|c}
             \!\!P_{\boldsymbol{x}_{k+1}\mid\boldsymbol{z}_k}(\cdot \mid z) &  
             \begin{aligned}
                 & \boldsymbol{x}_{k+1} = \bar{f}_k(z) + \boldsymbol{w}_k\\
                 &\mathbb{E}_{w\sim P_{\boldsymbol{w}_k}}[(w-\bar{\mu}_k)(w-\bar{\mu}_k)^\top] \!=\! \gamma_{2}\bar{\Sigma}_k
             \end{aligned}
        \end{array}\right\}.
    \end{equation*}
    iii) The objective function at $(\hat{p}^{*}_k, \rho^{*}_k)$ is
    \begin{equation*}
        -\frac{1}{2}\left[d_x\log(2\pi)+d_x+\log\operatorname{det}(\gamma_{2}\bar{\Sigma}_k)\right].
    \end{equation*}
\end{theorem}
\begin{proof}
Please see Appendix \ref{app:thm:one-step-DRPP-mp}.
\end{proof}
First, an immediate application of Theorem \ref{thm:one-step-DRPP-mp} is to utilize the explicit solution of $\hat{p}^{*}_k$ to develop a probabilistic prediction algorithm. Because this predictor solely focuses on the distributional robustness against the ambiguity of system noises, we name it Noise-DRPP. The pseudo-code is in Algorithm \ref{alg:noise-drpp}.

Second, notice that the optimal predictive pdf is unique, but the worst-case SDS is not. It contains any SDS that has a proper one-step state evolution and a noise whose first two moments satisfy an equation.

Finally, which is also our original goal, the objective function at $(\hat{p}^{*}_k, \rho^{*}_k)$ is an upper bound of the optimal objective of $\mathbf{P}_1$. Since the one-step upper bound does not depend on $z$, one can recursively use the Bellman equation \eqref{eq:bellman} to get an upper bound of the robust value function $V^{*}_k$.
\begin{theorem}[Upper bound of $V^{*}_k$]\label{cor:upper}
    Given an initial state-control pair $z\in\mathcal{Z}$ at time step $k\in\{0,\ldots,T-1\}$, an upper bound of the robust value function $V^{*}_k(z)$ is
    \begin{equation*}
        \begin{aligned}
            &V^{*}_k(z)\leq\sum_{t=k}^{T-1} -\frac{1}{2}\left[d_x\log(2\pi)+d_x+\log\operatorname{det}(\gamma_{2}\bar{\Sigma}_t)\right].
        \end{aligned}
    \end{equation*}
\end{theorem}
\begin{remark}
    We highlight that this upper bound can be computed offline because it only depends on the parameters of ambiguity sets $\mathcal{I}_{0:T-1}$. 
    % It can be used to verify whether the current ambiguity set is too conservative by comparing it with the real-time performance of a Noise-DRPP predictor in practice: if the scores always exceed this upper bound, one can suspect that the ambiguity set is too large. 
    Furthermore, the conservatism of this upper bound can be evaluated by the gap between it and another lower bound, as developed in the next section.
\end{remark}

\begin{algorithm}[h]
\renewcommand{\algorithmicrequire}{\textbf{Input:}}
\renewcommand{\algorithmicensure}{\textbf{Output:}}
\caption{Noise-DRPP $\mathscr{F}_{\text{Noise}}$} 
\label{alg:noise-drpp}
\begin{algorithmic}[1]
\Require time horizon $T$, ambiguity sets $\mathcal{I}_{0:T-1}$, control policy $\pi_{0:T-1}$, initial state $x_0$.
\State $s_0\leftarrow 0$.  \Comment{score at time step $0$}
\For {$k=0,\ldots,T-1$}
    \State Predictor updates ambiguity set $\mathcal{I}_k$ as \eqref{eq:ambiguity-set}.
    \State SDS generates control input $u_k \sim \pi_k(\cdot \mid x_k)$.
    \State $z_k \leftarrow (x_k, u_k)$, $\bar{\nu}_k \leftarrow \bar{f}_k(z_k)$.
    \State Predictor predicts $\hat{p}_k^{*} \sim \mathcal{N}\left(\bar{\nu}_k + \bar{\mu}_k, \;\gamma_{2}\bar{\Sigma}_k \right)$.
    \State SDS generates the next state $x_{k+1}$.
    \State $s_{k+1} \leftarrow s_k + \mathcal{L}(\hat{p}_k^{*}, x_{k+1})$.
\EndFor
\Ensure states $x_{0:T}$, predictions $\hat{p}^{*}_{0:T-1}$, scores $s_{0:T}$.  
\end{algorithmic} 
\end{algorithm}

\section{Value function Lower Bound}
To get a lower bound of $V^{*}_k(z)$, one can either shrink the feasible region of the outer maximization or enlarge the feasible region of the inner minimization. Because enlarging the ambiguity set does not essentially change the structure to alleviate the difficulty of a one-step DRPP $\mathbf{P}_2$, we choose to restrict $\hat{p}_k$ to a smaller pdf family by forcing the eigenvectors of the predictive covariance $\hat{\Sigma}_k$ to be the same as the nominal covariance $\bar{\Sigma}_k$. In this section, we elaborate on how the eigenvector restriction contributes to a well-performed algorithm and a lower bound.

\subsection{Problem Reformulation}
By adding auxiliary constraints to the predictor in \eqref{eq:qcqp}, one can still get a suboptimal prediction performance guarantee in the worst case, which is also a performance lower bound to the original one-step DRPP. Nevertheless, the choice of relaxation methods is delicate, which can greatly affect both the optimality gap and the computation efficiency.

Suppose the spectral decomposition for $\bar{\Sigma}_k$ is $Q_k \Lambda_k Q_k^\top$, where $Q_k=[\mathbf{v}_{1,k},\ldots,\mathbf{v}_{d_x,k}]$ is an orthogonal matrix whose columns are eigenvectors of $\bar{\Sigma}_k$ and $\Lambda_k = \operatorname{diag}\{\lambda_{1,k}, \ldots, \lambda_{d_x, k}\}$. Similarly, we let the spectral decomposition for $\hat{\Sigma}_k$ be $\hat{Q}_k \hat{\Lambda}_k \hat{Q}_k^\top$, where $\hat{Q}_k$ is an orthogonal matrix whose columns are eigenvectors of $\hat{\Sigma}_k$ and $\hat{\Lambda}_k = \operatorname{diag}\{\hat{\lambda}_{1,k}, \ldots, \hat{\lambda}_{d_x, k}\}$. Then, the eigenvector restriction for \eqref{eq:qcqp} is to impose the following constraint on the predictor:
\begin{equation}\label{eq:eigenvector-restriction}
    \hat{Q}_k = Q_k.
\end{equation}
To explain why we chose the eigenvector restriction, some supporting lemmas need to be derived. 
\begin{lemma}\label{lem:convex-max-bstar}
    The optimal $b^{*}$ to problem \eqref{eq:convex-max} is an eigenvector of $\frac{\alpha}{\|b^{*}\|_A} A + 2B$ with its $\ell^2$-norm being $\sqrt{\gamma_{0}(z)}$.
\end{lemma}
\begin{proof}
    % Please see Appendix \ref{app:lem:convex-max-bstar}.
    There are necessary optimality conditions for \eqref{eq:convex-max}:
    \begin{equation*}
        \left\{\begin{aligned}
            & \left(\frac{\alpha}{\|b\|_A} A + 2B -2sI\right)b = 0\\
            & \|b\|_2^2 = \gamma_{0}(z),
        \end{aligned}\right.
    \end{equation*}
    where the first one comes from the KKT condition and the second one holds because the objective increases monotonously with $\|b\|_2$.
    It follows that the optimal $b$ should be an eigenvector of $\frac{\alpha}{\|b\|_A} A + 2B$ with its $\ell^2$-norm being $\sqrt{\gamma_{0}(z)}$. 
\end{proof}
If the eigenvectors of $\frac{\alpha}{\|b^{*}\|_A} A + 2B$ have explicit expressions, one may have explicit expressions for $a^{*}, b^{*}$. Then, substituting these expressions into \eqref{eq:qcqp} can transform the original minimax optimization problem into a semidefinite-constrained maximization, where tractable optimization solvers may be available. However, it is always difficult to explicitly derive the eigenvectors when $\hat{Q}_k$ does not identify with $Q_k$. 
% Moreover, there is
% \begin{equation*}
%     \begin{aligned}
%         &0 = b^\top\left(\frac{\alpha}{\|b\|_A} A + 2B -2sI\right)b= \alpha\|b\|_A + 2\|b\|_B^2 - 2s\|b\|^2\\
%         \Rightarrow & \|b\|_B^2 = s\gamma_0 - \frac{\alpha}{2}\|b\|_A.
%     \end{aligned}
% \end{equation*}

\begin{lemma}\label{lem:convex-max-relax}
    When $\hat{Q}_k = Q_k$ is supplemented to the constraints of \eqref{eq:convex-max}, the solution is
    \begin{equation*}
        b^{*} = \sqrt{\gamma_{0}(z)}\mathbf{v}_{j_k,k},
    \end{equation*}
    where $j_k=\arg\max\limits_{i}(2\sqrt{\gamma_{0}(z)\gamma_1\lambda_{i,k}}+\gamma_{0}(z))\hat{\lambda}_{i,k}^{-1}$.
\end{lemma}
\begin{proof}
    % Please see appendix \ref{app:lem:convex-max-relax}.
    When $\hat{Q}_k = Q_k$ is satisfied, there is
    \begin{equation*}
        \frac{\alpha}{\|b\|_A} A + 2B =  Q_k\left[\frac{\alpha}{\|b\|_A}\hat{\Lambda}_k^{-2}\Lambda_k + \hat{\Lambda}_k^{-1}\right]Q_k^\top,
    \end{equation*}
    whose eigenvectors are the columns of $Q_k$. Lemma \ref{lem:convex-max-bstar} indicates that the optimal $b$ can be expressed as $\sqrt{\gamma_{0}(z)}\mathbf{v}_i$ for certain $i\in\{1,\ldots, d_x\}$, then we have the objective of \eqref{eq:convex-max} is $(2\sqrt{\gamma_{0}(z)\gamma_1\lambda_{i,k}}+\gamma_{0}(z))\hat{\lambda}_{i,k}^{-1}$. Since $j_k=\arg\max_{i}(2\sqrt{\gamma_{0}(z)\gamma_1\lambda_{i,k}}+\gamma_{0}(z))\hat{\lambda}_{i,k}^{-1}$, problem \eqref{eq:convex-max} is solved as $b^{*} = \sqrt{\gamma_{0}(z)}\mathbf{v}_{j_k}$.
\end{proof}

Utilizing the expression of $b^{*}$ in Lemma \ref{lem:convex-max-relax}, the minimax optimization $\eqref{eq:qcqp}$ under eigenvector restriction \eqref{eq:eigenvector-restriction} can be transformed into
\begin{equation*}\label{eq:final-convex}
    \begin{aligned}
        (\mathbf{P}_3):&\min_{\hat{\Lambda}_k\succeq 0, j_k} \sum_{i=1}^{d_x} \left[-\log(\hat{\lambda}_{i,k}^{-1}) +\gamma_{2}\lambda_{i,k}\hat{\lambda}_{i,k}^{-1}\right]+\\
        &\qquad \left(2\sqrt{\gamma_{0}(z)\gamma_1\lambda_{j_k,k}}+\gamma_{0}(z)\right)\hat{\lambda}_{j_k,k}^{-1}\\
        &\text{s.t. } \hat{\lambda}_{j_k,k} \leq \frac{2\sqrt{\gamma_{0}(z)\gamma_1\lambda_{j_k,k}}\!+\! \gamma_{0}(z)}{2\sqrt{\gamma_{0}(z)\gamma_1\lambda_{i,k}}\!+\! \gamma_{0}(z)} \hat{\lambda}_{i,k} \\
        &\qquad\qquad\qquad\qquad\qquad\qquad\qquad \forall i\in\{1,\ldots,d_x\}.
    \end{aligned}
\end{equation*}

\subsection{Eig-DRPP}
Notice that optimization $\mathbf{P}_3$ is convex for any fixed $j_k$. Therefore, $\mathbf{P}_3$ can be solved by taking the minimum of $d_x$ different convex optimization problems, which is numerically efficient.
Let $\left(\hat{\Lambda}_k^{*}=\operatorname{diag}\{\hat{\lambda}_{1,k}^{*}, \ldots, \hat{\lambda}_{d_x,k}^{*}\}, j_k^{*}\right)$ be the solution of $\mathbf{P}_3$, we summarize the solution of a one-step DRPP under eigenvector restriction as follows.
\begin{theorem}[Eig-DRPP]\label{thm:one-step-DRPP-relax}
    Given $\boldsymbol{z}_k = z \in \mathcal{Z}$ at time step $k\in\{0,\ldots,T-1\}$, if the probabilistic predictor is constrained by the eigenvector restriction \eqref{eq:eigenvector-restriction}, the solution to the one-step DRPP $\mathbf{P}_1$ is:\\
    i) The optimal predictive pdf $\hat{p}^{*}_k$ is
    \begin{equation*}
        \hat{p}_k^{*} \sim \mathcal{N}\left(\bar{\nu}_k + \bar{\mu}_k, \;Q_k\hat{\Lambda}_k^{*} Q_k^\top \right).
    \end{equation*}
    ii) The worst-case conditional measure $\rho_k^*$ belongs to the set
    \begin{equation*}
        \left\{\!\begin{array}{c|c}
             \!\!P_{\boldsymbol{x}_{k+1}\mid\boldsymbol{z}_k}(\cdot \mid z)&  
             \begin{aligned}
                 &\boldsymbol{x}_{k+1} = f_k(z) + \boldsymbol{w}_k\\
                 & f_k(z) = \sqrt{\gamma_{0}(z)}\mathbf{v}_{j_k^{*}\!,k} + \bar{\nu}_k\\
                 &\mathbb{E}[\boldsymbol{w}_k] = \bar{\mu}_k+\sqrt{\gamma_1\lambda_{j_k^{*}\!, k}}\mathbf{v}_{j_k^{*}\!,k}\\
                 &\operatorname{Cov}[\boldsymbol{w}_k] = \gamma_{2}\bar{\Sigma}_k-\gamma_1\lambda_{j_k^{*}\!, k} \mathbf{v}_{j^{*}\!, k}\mathbf{v}_{j_k^{*}\!, k}^\top
             \end{aligned}
        \end{array}\right\}.
    \end{equation*}
    iii) The objective function at $(\hat{p}^{*}_k, \rho^{*}_k)$ is
    \begin{equation*}
    \begin{aligned}
       &-\frac{1}{2}\left\{d_x\log(2\pi)+\sum_{i=1}^{d_x} \left[\log(\hat{\lambda}_{i,k}^{*}) +\gamma_{2}\frac{\lambda_{i,k}}{\hat{\lambda}_{i,k}^{*}}\right]\right.\\
       & \left.\qquad+\frac{2\sqrt{\gamma_{0}(z)\gamma_1\lambda_{j_k^{*}\!,k}}+\gamma_{0}(z)}{\hat{\lambda}_{j_k^{*}\!,k}^{*}}\right\}.
    \end{aligned}
    \end{equation*}
\end{theorem}
\begin{proof}
    % Please see Appendix \ref{app:thm:one-step-DRPP-relax}.
    The proof is completed by combining the results in Lemma \ref{lem:primal-solution-1}, Lemma \ref{lem:QCLP}, and the solutions of $\mathbf{P}_3$.
\end{proof}

First, an immediate application of Theorem \ref{thm:one-step-DRPP-relax} is to utilize the explicit solution of $\hat{p}^{*}_k$ to develop a probabilistic prediction algorithm, named Eig-DRPP. The pseudo-code is presented in Algorithm \ref{alg:eig-drpp}. 
Second, although the optimal predictive pdf is unique, the worst-case conditional measure is not. Compared to Theorem \ref{thm:one-step-DRPP-mp}, there are relatively more restrictions on the set of worst-case conditional measures.
Third, when $\gamma_{0}=0$, i.e., the predictor has no ambiguity of $f_k$, the solution of $\mathbf{P}_3$ is $\hat{\lambda}_{i,k}^{*} = \gamma_{2}\lambda_{i,k}$ and $j_k^{*}$ can be any feasible index. In this case, both the optimal predictive pdf and objective value in Theorem \ref{thm:one-step-DRPP-relax} coincide with those in Theorem \eqref{thm:one-step-DRPP-mp}. However, the sets of worst-case conditional measures are not the same, due to the extra constraint for the expectation in Theorem \ref{thm:one-step-DRPP-relax}. 

Finally, the objective function at $(\hat{p}^{*}_k, \rho^{*}_k)$ is a lower bound of the optimal objective of $\mathbf{P}_1$. Since the objective value depends on $z$ through $\gamma_{0}(z)$, one cannot directly use the Bellman equation \eqref{eq:bellman} recursively to get a lower bound of $V^{*}_k$. 
\begin{theorem}[Lower bound of $V^{*}_k$ and optimality gap of $\mathscr{F}_{\text{Eig}}$]\label{cor:lower}
    Given an initial state-control pair $z\in\mathcal{Z}$ at time step $k\in\{0,\ldots,T-1\}$, and an upper bound of $\gamma_{0}$ such that $\gamma_{0}(z) \leq \Gamma_{0}\!\in\!\mathbb{R}_{+}, \forall z\in\mathcal{Z}$, there is\\
    i) A lower bound of the robust value function $V_k^*(z)$ is
    \begin{equation*}
        \begin{aligned}
            V^{*}_k(z)\geq & \sum_{t=k}^{T-1} \!-\frac{1}{2}\!\left\{d_x\log(2\pi)\!+\!\sum_{i=1}^{d_x} \left[\log(\hat{\lambda}_{i,t}^{*}) \!+\!\gamma_{2}\frac{\lambda_{i,t}}{\hat{\lambda}_{i,t}^{*}}\right]\right.\\
           & \left.\qquad+\frac{2\sqrt{\Gamma_{0}\gamma_1\lambda_{j_t^{*}\!,t}}+\Gamma_{0}}{\hat{\lambda}_{j_t^{*}\!,t}^{*}}\right\},
        \end{aligned}
    \end{equation*}
    where $\hat{\lambda}_{1,t}^{*}, \ldots, \hat{\lambda}_{d_x,t}^{*}, j_t^{*}$ is the solution to $\mathbf{P}_3$ with $k$ being replaced by $k$ and all $\gamma_{0}(z)$ being replaced by $\Gamma_{0}$.\\
    ii) The optimality gap between Eig-DRPP $\mathscr{F}_{\text{Eig}}$ and the global optimal predictor is upper bounded as follows,
    \begin{equation}\label{eq:optimality_gap}
        \begin{aligned}
            &V_k^{*}(z) - V_k^{\mathscr{F}_{\text{Eig}}}(z)  \leq \sum_{t=k}^{T-1} -\frac{1}{2}\bigg\{\underbrace{d_x + \log\operatorname{det}(\gamma_{2}\bar{\Sigma}_t)}_{\circled{1}}\\
             & \underbrace{-\sum_{i=1}^{d_x} \left[\log(\hat{\lambda}_{i,t}^{*}) \!+\!\gamma_{2}\frac{\lambda_{i,t}}{\hat{\lambda}_{i,t}^{*}}\right]
          -\frac{2\sqrt{\Gamma_{0}\gamma_1\lambda_{j_t^{*}\!,t}}+\Gamma_{0}}{\hat{\lambda}_{j_t^{*}\!,t}^{*}}}_{\circled{2}}\bigg\}.
        \end{aligned}
    \end{equation}
\end{theorem}
\begin{proof}
    i) Notice that as $\gamma_{0}(z)$ increases, the ambiguity set enlarges, thus the optimal value of one-step DRPP under eigenvector constraint should monotonously decrease. Since the upper bound of $\gamma_{0}(z)$, i.e., $\Gamma_{0}$, is independent of $z$, one can still use \eqref{eq:bellman} to derive a lower bound.
    
    ii) Subtracting this lower bound from the upper bound in Theorem \ref{cor:upper}, one immediately gets an upper bound of the optimality gap for Eig-DRPP. As illustrated in \eqref{eq:optimality_gap}, $\circled{1}$ comes from the upper bound, which is determined by $\gamma_{2}$ and the determinant of $\bar{\Sigma}_t$. While $\circled{2}$ comes from the lower bound, which is jointly determined by $\gamma_{0}, \gamma_1, \gamma_{2}$ and the eigenvalues of $\bar{\Sigma}_t$ by solving $\mathbf{P}_3$.
\end{proof}
% \begin{remark}
%     Similar to the upper bound, this lower bound can also be computed offline because it only depends on the ambiguity sets $\mathcal{I}_{0:T-1}$. It can be used to verify whether the current ambiguity set is too optimistic by comparing it with the real-time performance of an Eig-DRPP predictor in practice: if the scores always stay lower than this lower bound, one can suspect that the ambiguity set is too small.
% \end{remark}
The proposed DRPP predictors can be naturally integrated into SMPC. In particular, the Noise-DRPP predictor provides a closed-form Gaussian predictive model that can be used to construct distributionally robust chance constraints. The Eig-DRPP predictor, while computationally more involved, can serve as an offline surrogate to characterize the trade-off between robustness and conservatism in SMPC design.
\begin{algorithm}[h]
\renewcommand{\algorithmicrequire}{\textbf{Input:}}
\renewcommand{\algorithmicensure}{\textbf{Output:}}
\caption{Eig-DRPP $\mathscr{F}_{\text{Eig}}$} 
\label{alg:eig-drpp}
\begin{algorithmic}[1]
\Require time horizon $T$, ambiguity sets $\mathcal{I}_{0:T-1}$, control policy $\pi_{0:T-1}$, initial state $x_0$.
\State $s_0\leftarrow 0$.  \Comment{score at time step $0$}
\For {$k=0,\ldots,T-1$}
    \State Predictor updates ambiguity set $\mathcal{I}_k$ as \eqref{eq:ambiguity-set}.
    \State Do spectral decomposition for $\bar{\Sigma}_k = Q_k \Lambda_k Q_k^\top$ .
    \State SDS generates control input $u_k \sim \pi_k(\cdot \mid x_k)$.
    \State $z_k \leftarrow (x_k, u_k)$, $\bar{\nu}_k \leftarrow \bar{f}_k(z_k)$.
    \State $j_k^{*} \leftarrow 0$, $\operatorname{val}^{*} \leftarrow \infty$, $\hat{\Lambda}_k^{*} \leftarrow \Lambda_k$.
    \For {$j=1,2,\ldots,d_x$}  
        \State Fix $j_k = j$, solve $\mathbf{P}_3$, where the optimizer is $\hat{\Lambda}_k^{(j)}$ and the optimal value is $\operatorname{val}^{(j)}$. \Comment{For each $j$, this step is a convex optimization}
        \If {$\operatorname{val}^{(j)}<\operatorname{val}^{*}$}
            \State $j_k^{*} \leftarrow j$, $\operatorname{val}^{*} \leftarrow \operatorname{val}^{(j)}$, $\hat{\Lambda}_k^{*} \leftarrow \hat{\Lambda}_k^{(j)}$.
        \EndIf
    \EndFor
    \State Predictor predicts $\hat{p}_k^{*} \sim \mathcal{N}\left(\bar{\nu}_k + \bar{\mu}_k, \;Q_k\hat{\Lambda}_k^{*} Q_k^\top \right)$.
    \State SDS generates the next state $x_{k+1}$.
    \State $s_{k+1} \leftarrow s_k + \mathcal{L}(\hat{p}_k^{*}, x_{k+1})$.
\EndFor
\Ensure states $x_{0:T}$, predictions $\hat{p}^{*}_{0:T-1}$, scores $s_{0:T}$.  
\end{algorithmic} 
\end{algorithm}
\section{Experiments}
In this section, a series of experiments is conducted to explore three questions. i) Given an SDS subject to an ambiguity set, what are the performance advantages of Noise-DRPP and Eig-DRPP compared to a nominal predictor? ii) How much will the prediction performances be influenced by control strategies? iii) How can DRPP predictors be applied to providing high probability confidence regions of future states?
\begin{figure*}[h]
    \centering
    \subfigure[SDS \eqref{eq:LTI} under zero control]{
        \includegraphics[width = 0.31\textwidth]{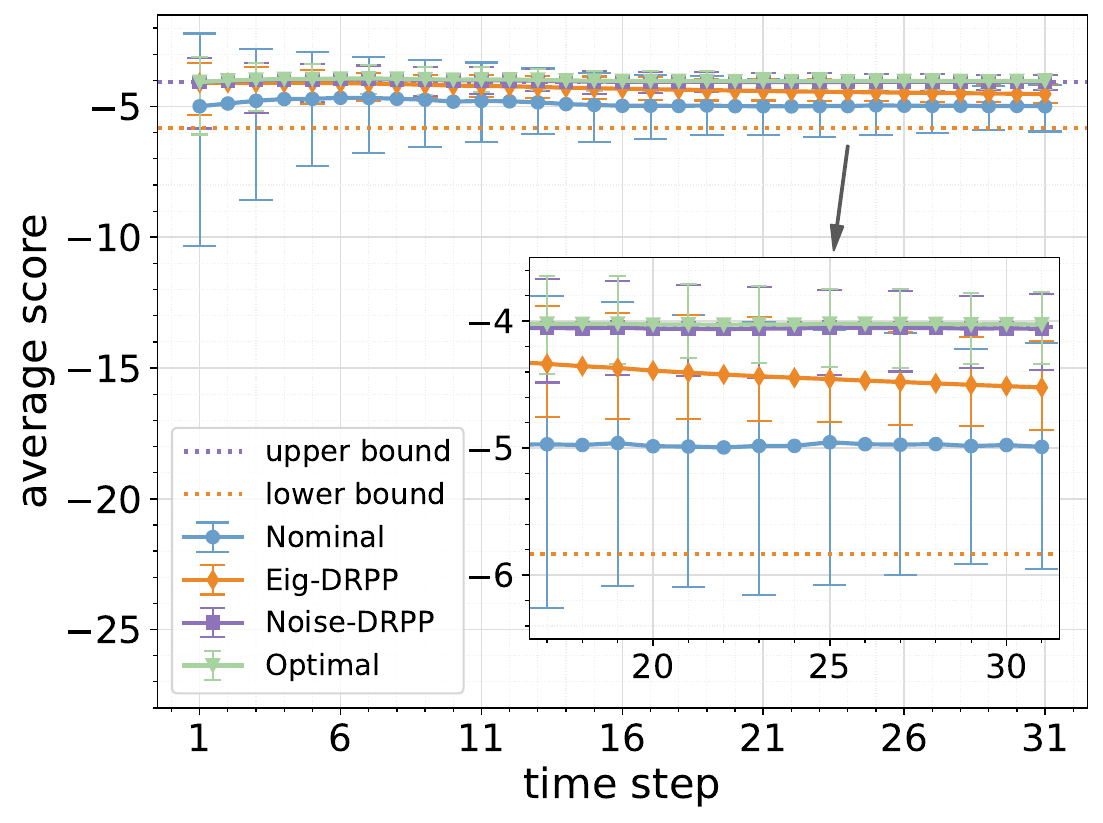}
        \label{auto:1:eb}
    }
    \subfigure[SDS \eqref{eq:LTV} under zero control]{
        \includegraphics[width = 0.31\textwidth]{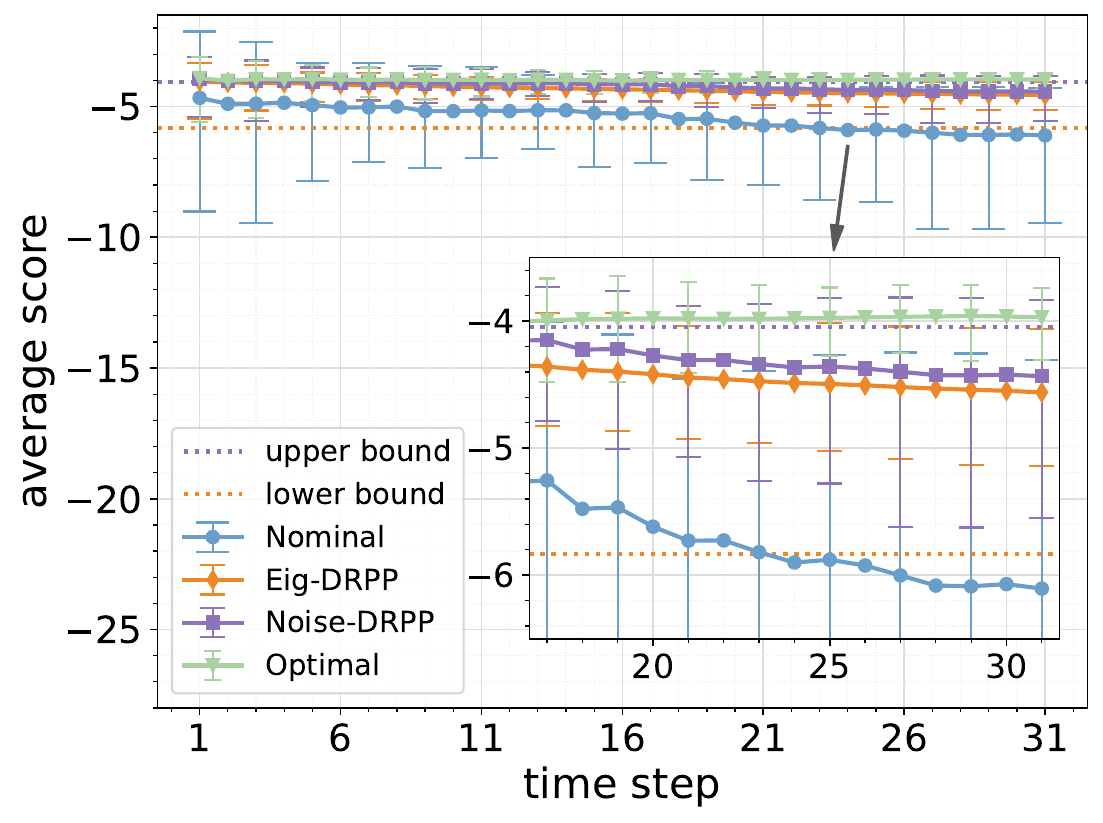}
        \label{auto:2:eb}
    }
    \subfigure[Adversarial SDS under zero control]{
        \includegraphics[width = 0.31\textwidth]{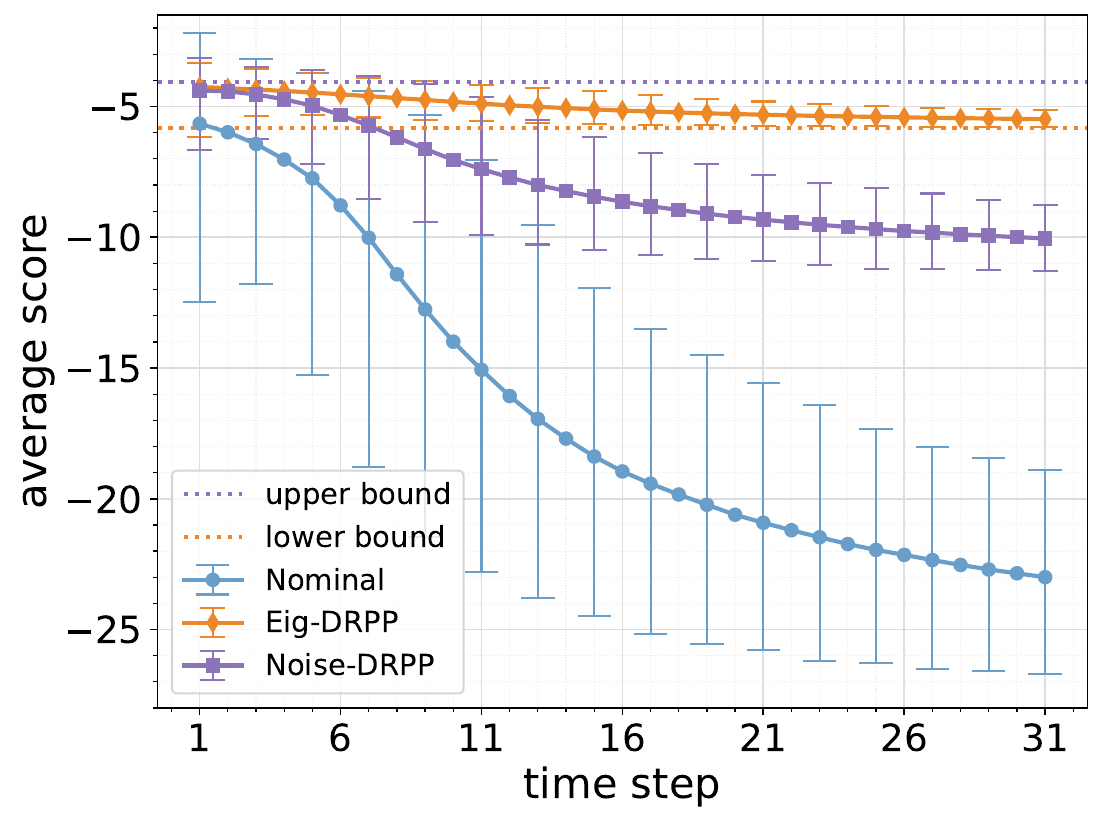}
        \label{auto:3:eb}
    }
    \subfigure[SDS \eqref{eq:LTI} under LQR control]{
        \includegraphics[width = 0.31\textwidth]{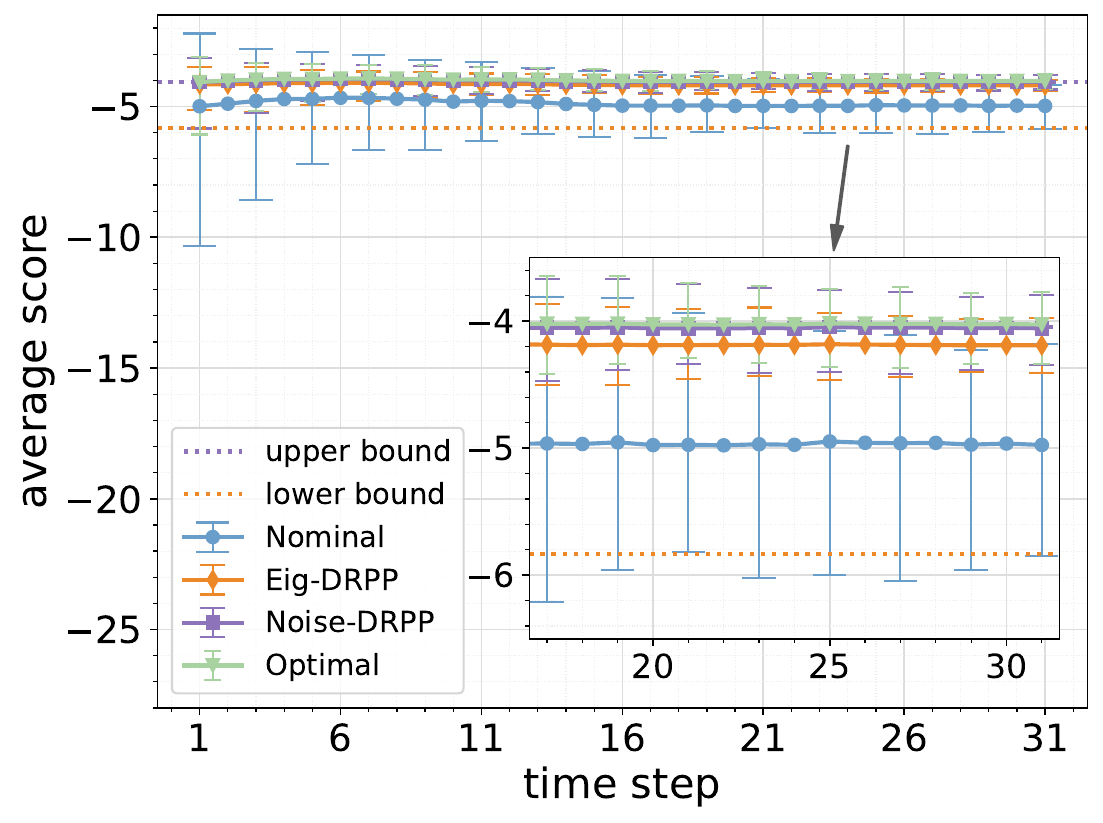}
        \label{regulation:1:eb}
    }
    \subfigure[SDS \eqref{eq:LTV} under LQR control]{
        \includegraphics[width = 0.31\textwidth]{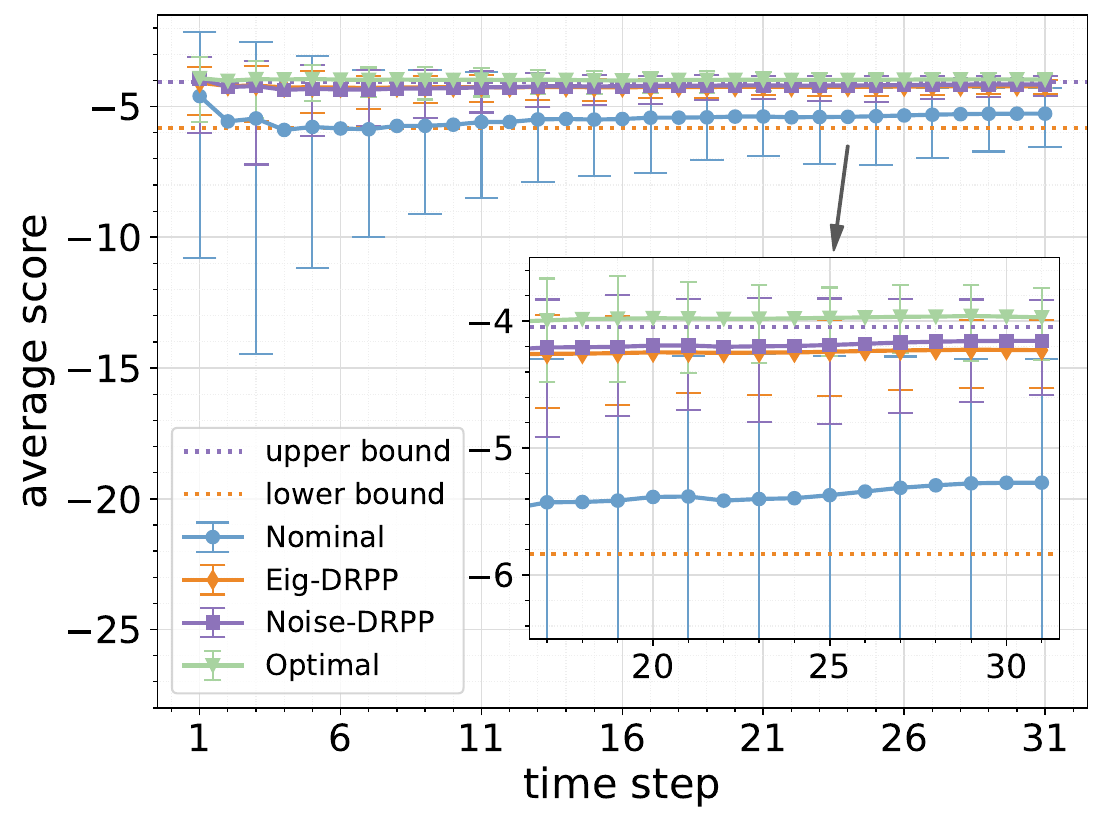}
        \label{regulation:2:eb}
    }
    \subfigure[Adversarial SDS under LQR control]{
        \includegraphics[width = 0.31\textwidth]{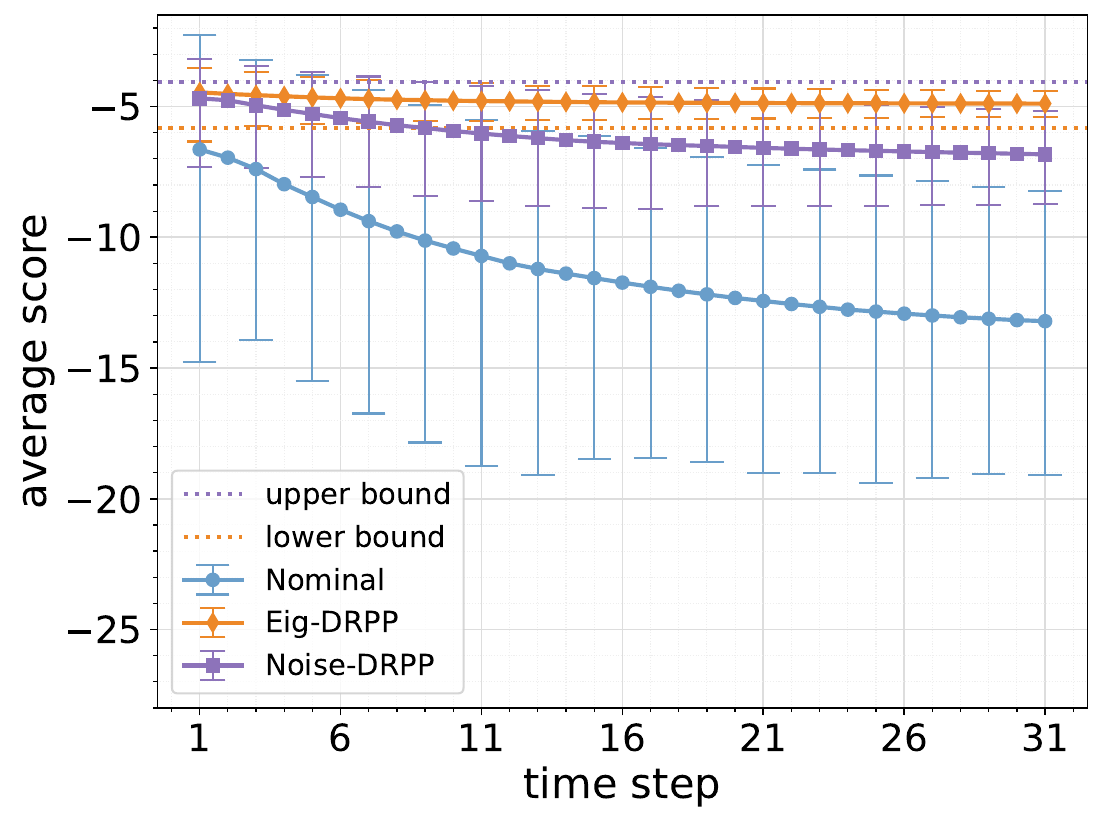}
        \label{regulation:3:eb}
    }
    \caption{Prediction performance of different probabilistic predictors on different SDSs under different control strategies. }
    \label{fig:performance}
\end{figure*}
\subsection{Experiment Setup}
\paragraph*{Ambiguity set} Conditioned on a state-control pair $z_k$ at time step $k\in\{0, 1, \cdots, 31\}$, the nominal state evolution function is given as
\begin{equation*}
    \bar{f}_k(\boldsymbol{x}_k, \boldsymbol{u}_k) = \begin{pmatrix}
    1 & 0.1\\
    0 & 1\\
\end{pmatrix} \boldsymbol{x}_k + \begin{pmatrix}
    1 & 0\\ 
    0 & 1
\end{pmatrix} \boldsymbol{u}_k,
\end{equation*}
and the uncertainty between $\bar{f}_k$ and the true $f_k$ is quantified by $\gamma_{0}(z_k) = \min\{0.3\|z_k\|_2, 5\}^2$.
The nominal mean and covariance of the noise $\boldsymbol{w}_k$ conditioned on $z_k$ are $\bar{\mu}_k = \begin{pmatrix}
        0 \\
        0
    \end{pmatrix}, 
    \bar{\Sigma}_k = \begin{pmatrix}
        1 & 0.5 \\
        0.5 & 1.5
    \end{pmatrix}$, and the uncertainty between $\bar{\mu}_k, \bar{\Sigma}_k$ and the real $\mu_k, \Sigma_k$ are quantified by $\gamma_1=0.5$ and $\gamma_{2}=3$, respectively. Given the ambiguity set specified above, the real SDS belongs to $\mathcal{I}_k$. Although the nominal model is linear time-invariant (LTI), the real model does not have to be LTI. 

\paragraph*{Simulation mechanism} Three different simulation mechanisms are conducted to generate the underlying SDSs.

i) Let $\boldsymbol{\alpha}_1, \boldsymbol{\alpha}_2$ be two independent random variables subject to the uniform distribution on $[-1,1]$, and let $\boldsymbol{\alpha}_3$ be another independent random variable subject to the uniform distribution on $[-1,1]^2$. We randomly generates $\alpha_i \sim p_{\boldsymbol{\alpha}_i}$ for $i=1,2,3$, then we simulate the SDS as
\begin{equation}\label{eq:LTI}
    \begin{aligned}
        \boldsymbol{x}_{k+1} = \begin{pmatrix}
    1 & 0.1+0.3\alpha_1\\
    0 & 1\\
\end{pmatrix} \boldsymbol{x}_k + \begin{pmatrix}
    1 & 0.3\alpha_2\\ 
    0 & 1
\end{pmatrix} \boldsymbol{u}_k + \boldsymbol{w}_k,
    \end{aligned}
\end{equation}
where $\boldsymbol{w}_k \sim \mathcal{N}\left(0.5\alpha_3, 3\bar{\Sigma}_k - 0.25\alpha_3\alpha_3^\top\right)$.

ii) For each $i=1,2,3$, let $\{\boldsymbol{\alpha}_{i,k}\}_{k=0}^{31}$ be random variables that are independent and identically distributed to $p_{\boldsymbol{\alpha}_i}$. We randomly generates $\alpha_{i,k}\sim p_{\boldsymbol{\alpha}_{i,k}}$ for each $i$ and $k$, then we simulate the SDS as 
\begin{equation}\label{eq:LTV}
    \begin{aligned}
        \boldsymbol{x}_{k+1} = \begin{pmatrix}
    1 & 0.1+0.3\alpha_{1,k}\\
    0 & 1\\
\end{pmatrix} \boldsymbol{x}_k + \begin{pmatrix}
    1 & 0.3\alpha_{2,k}\\ 
    0 & 1
\end{pmatrix} \boldsymbol{u}_k + \boldsymbol{w}_{k},
    \end{aligned}
\end{equation}
where $\boldsymbol{w}_k \sim \mathcal{N}(0.5\alpha_{3,k}, 3\bar{\Sigma}_k - 0.25\alpha_{3,k}\alpha_{3,k}^\top)$.

iii) (Adversarial) In this case, the underlying SDS is allowed to adversarially choose an SDS in the ambiguity set at each step to degrade the prediction performance. Specifically, at each step, after a predictive pdf has been output by the predictor, the adversarial SDS uses the predicted covariance to choose a worst-case SDS as described in Theorem \ref{thm:one-step-DRPP-relax}.

\paragraph*{Predictors comparison} Two different control strategies are considered, the zero input and the linear quadratic regulation (LQR), where the zero input strategy means no control input is imposed, and the LQR control strategy is defined as a standard linear quadratic regulator based on the nominal linear model with $Q$ and $R$ all set as identity matrices. 
Under each control strategy and SDS scenario, we randomly generate 1,000 trajectories starting from the same initial state $x_0 = [2,1]$. 
% For the LTI and LTV mechanisms, since the ground-truth dynamics are known before prediction, we have an optimal predictor as the least upper bound of prediction performance. 
At each time step $k$, let the nominal predictor predicts $\bar{p}_k \sim \mathcal{N}(\bar{\nu}_k+\bar{\mu}_k, \bar{\Sigma}_k)$, and the optimal predictor predicts $p_k \sim \mathcal{N}(\nu_k+\mu_k, \Sigma_k)$.
Then, we compare the prediction performance among the nominal predictor, Noise-DRPP, Eig-DRPP, and the optimal predictor. 
Notice that for the adversarial mechanism, there is no explicitly predefined optimal predictor because the real SDS depends on the predictor. Therefore, we should compare the prediction performance among the nominal predictor, Noise-DRPP, and Eig-DRPP. 
% In fact, as guaranteed by Theorem \ref{thm:one-step-DRPP-relax}, when predictors are constrained by the eigenvector restriction, Eig-DRPP is the optimal predictor under this adversarial scenario.
\begin{figure*}[t]
    \centering
    \subfigure[SDS \eqref{eq:LTI} under zero control]{
        \includegraphics[width = 0.31\textwidth]{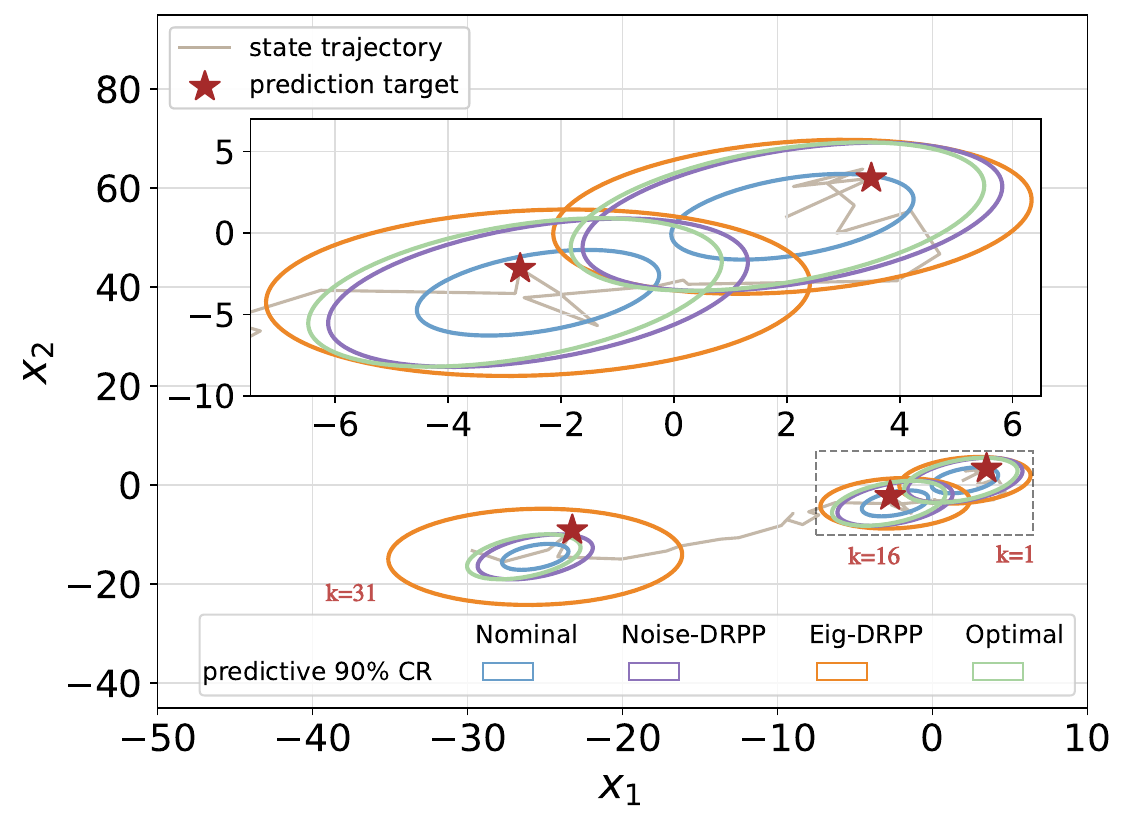}
        \label{auto:1:traj}
    }
    \subfigure[SDS \eqref{eq:LTV} under zero control]{
        \includegraphics[width = 0.31\textwidth]{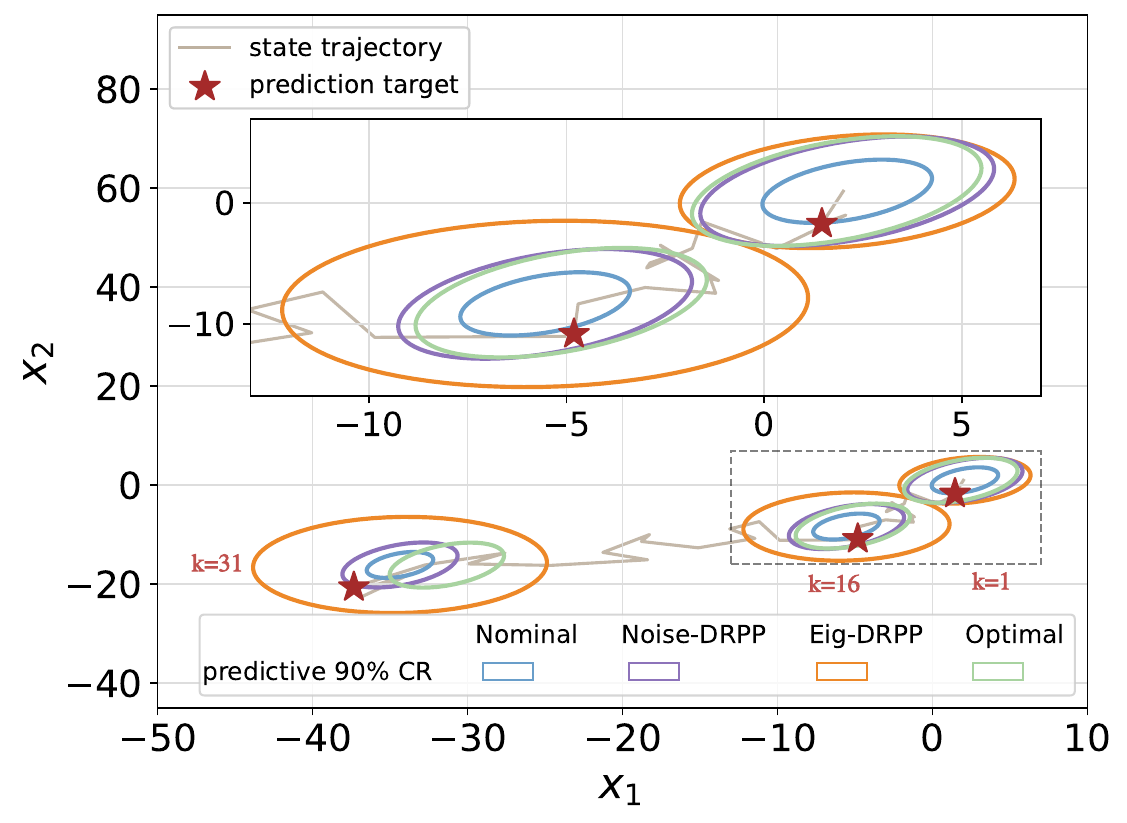}
        \label{auto:2:traj}
    }
    \subfigure[Adversarial SDS under zero control]{
        \includegraphics[width = 0.305\textwidth]{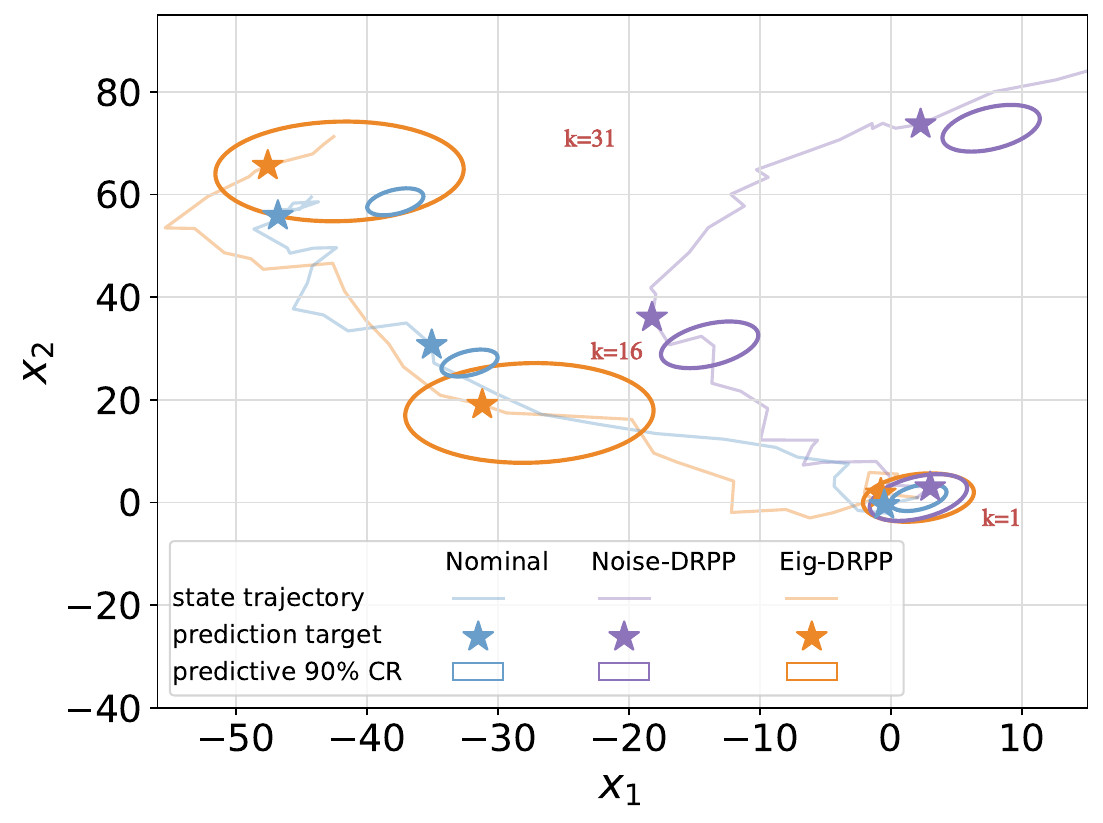}
        \label{auto:3:traj}
    }
    \subfigure[SDS \eqref{eq:LTI} under LQR control]{
        \includegraphics[width = 0.31\textwidth]{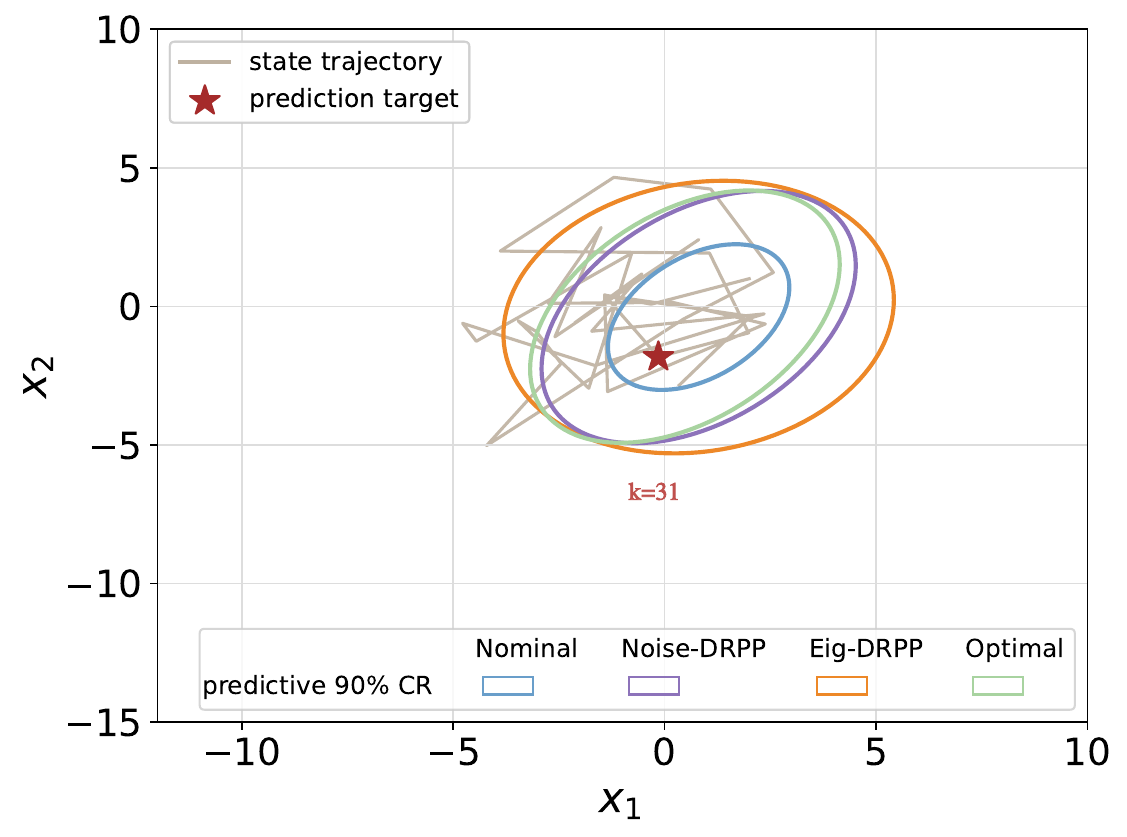}
        \label{regulation:1:traj}
    }
    \subfigure[SDS \eqref{eq:LTV} under LQR control]{
        \includegraphics[width = 0.31\textwidth]{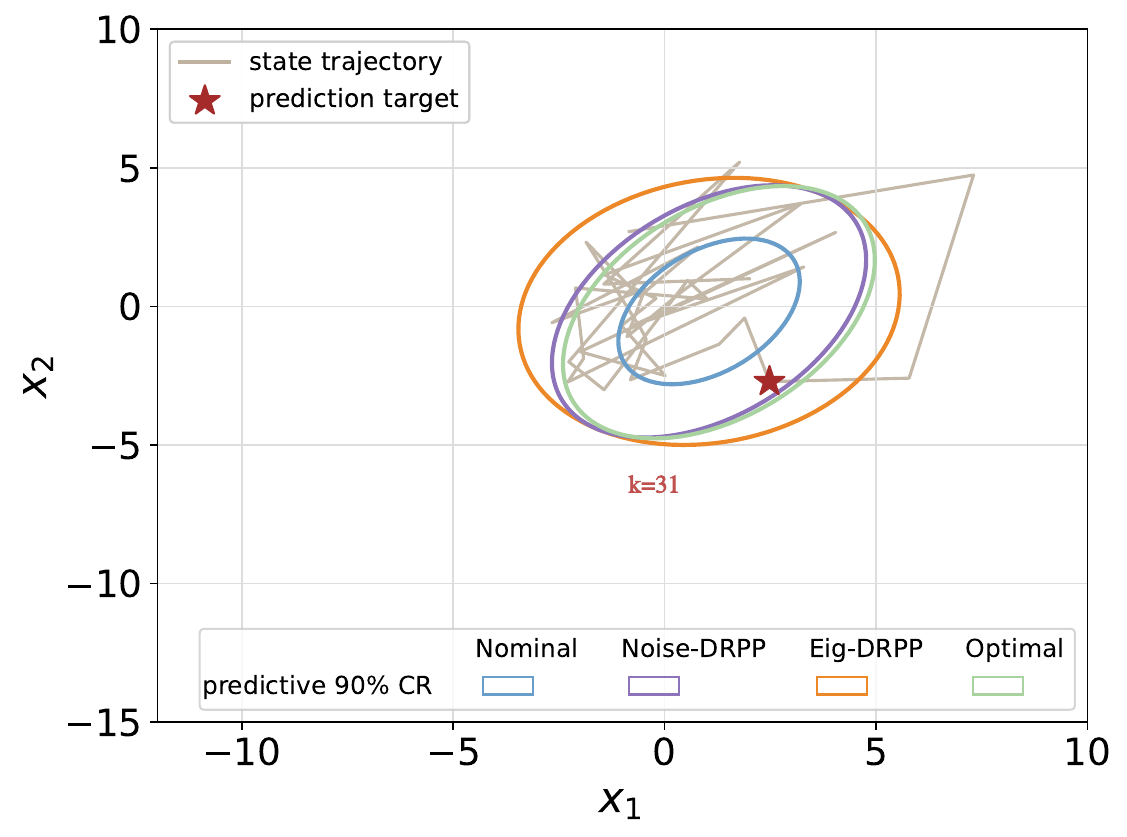}
        \label{regulation:2:traj}
    }
    \subfigure[Adversarial SDS under LQR control]{
        \includegraphics[width = 0.31\textwidth]{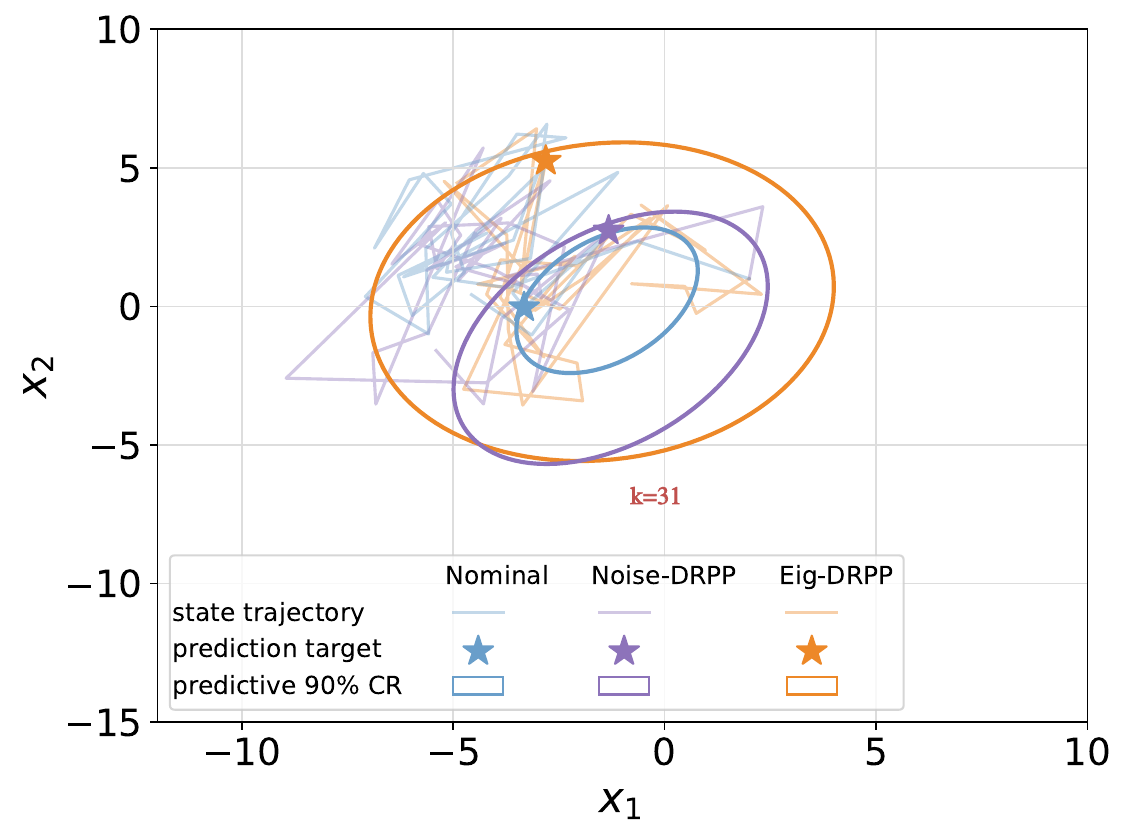}
        \label{regulation:3:traj}
    }
    \caption{Predictive 90\% confidence regions of probabilistic predictors for different SDSs under different control strategies.}
    \label{fig:cr}
\end{figure*}
\subsection{Results and Discussions}
In Fig. \ref{fig:performance}, the prediction performances of different probabilistic predictors on different SDSs under different control strategies are visualized and compared with each other. For each SDS, we evaluate the temporal average scores of each predictor at each time step on 1,000 randomly sampled trajectories. At each step, the mean of 1,000 average scores is plotted as a dot, which is an approximation to the expected average score. The 5\% and 95\% percentiles of these average scores are plotted as an error bar, reflecting how concentrated the average scores are.

\paragraph*{Performance advantages} Under each setting of SDS and control strategy, both Noise-DRPP and Eig-DRPP possess larger expected average scores at each time step compared to the nominal predictor. Since the ambiguity set for our simulation does not have too large an uncertainty, the scores of both DRPP predictors are close to the optimal predictor. Additionally, both DRPP predictors have more concentrated average scores at each time step. As the uncertainty of the real underlying SDS increases from LTI to adversarial, the performance gap between the DRPP predictors and the nominal predictor increases. 

% \paragraph*{Upper and lower bounds}
% When the ambiguity set is not assured to be accurate, we can use the theoretical upper and lower bounds to adjust the ambiguity set, and then improve the performance of DRPP predictors. We provide two simple but efficient rules to adjust the ambiguity set in practice:
% i) If the expected score of any DRPP predictor is always larger than the upper bound, it means that the true SDS is less adversarial than the worst-case SDS in $\mathcal{I}_k$, thus the predictor can use a smaller ambiguity set for better scores. 
% ii) If the expected score of an Eig-DRPP predictor is always smaller than the lower bound, it means that the true SDS is more adversarial than the worst-case SDS in $\mathcal{I}_k$, thus the Eig-DRPP predictor can use a larger ambiguity set for better scores.

\paragraph*{Influence of control strategies}
When a predictor has no ambiguity about the system's state evolution function, predicting the states is equivalent to predicting the noises. In this case, the choice of control strategies does not affect any probabilistic predictor. 
However, when a predictor does have ambiguity of the system's state evolution function, different control strategies can greatly influence the performance of a probabilistic predictor. Since the ambiguity set at each step directly depends on the state-input pair, different control strategies will lead to quite different state-input trajectories, thus leading to quite different prediction performances. We can observe this phenomenon by comparing (a-c) with (d-f) in Fig. \ref{fig:performance} respectively. An intuitive explanation is that the LQR control strategy regulates the system states towards the neighborhoods of $0$, where the ambiguity of one-step state evolution is smaller than other states, thus the DRPP predictors perform better.

\subsection{Application: Predict Robust Confidence Regions}
An immediate application of probabilistic prediction is to predict high-probability confidence regions for future states. Since the predictive pdfs of both Noise-DRPP and Eig-DRPP are Gaussian, one can easily derive an elliptical $\beta$-confidence region for any $\beta\in(0,1)$. In Fig. \ref{fig:cr}, we have visualized the predictive 90\% confidence regions for each predictor at certain time steps of a randomly generated trajectory. For the zero control scenario, we visualize the predictions for time steps $k=1,16,31$, while for the LQR control scenario, we visualize the prediction for time step $k=31$. 

Visualizing the behaviors of different probabilistic predictors on a trajectory, one can intuitively explain why DRPP predictors outperform the nominal predictor.
The nominal predictor is overly confident, where the predictive pdfs are always too sharp to enclose the future states in their 90\% CRs. In contrast, DRPP predictors have taken the ambiguity of noises into account, whose CRs are more robust to contain the prediction targets most of the time. If we further compare Noise-DRPP with Eig-DRPP, it seems that Eig-DRPP is more intelligent in adaptively changing the ellipses' shapes. This is because Eig-DRPP has also considered the ambiguity of one-step state evolution, and the eigenvalues of predictive covariances are attained based on minimax optimization.

\subsection{Extension: Nonlinear System}
We consider a 6-degree-of-freedom serial manipulator with revolute joints using the standard Euler-Lagrange rigid-body equations \cite{lynchModernRoboticsMechanics2017}. The coordinates are $q \in \mathbb{R}^6$ and velocities $\dot{q} \in \mathbb{R}^6$. The continuous-time dynamics are
\begin{equation*}
    M(q) \ddot{q}+C(q, \dot{q}) \dot{q}+g(q)=\tau,
\end{equation*}
where $M(q)$ is the inertia matrix, $C(q, \dot{q})\dot{q}$ represents the Coriolis forces, $g(q)$ is teh gravity vector, and $\tau \in \mathbb{R}^6$ denotes applied joint torques. The mass $m_i = 0.8$kg, link length $l_i = 0.5$m, and link center of mass (COM) distance $l_{c, i} = 0.25$m for each link $i=1, \ldots, 6$. 
We write the full state $x=\left[q^{\top}, \dot{q}^{\top}\right]^{\top} \in \mathbb{R}^{12}$.

\begin{figure}[htb]
    \centering
    \includegraphics[width=0.4\textwidth]{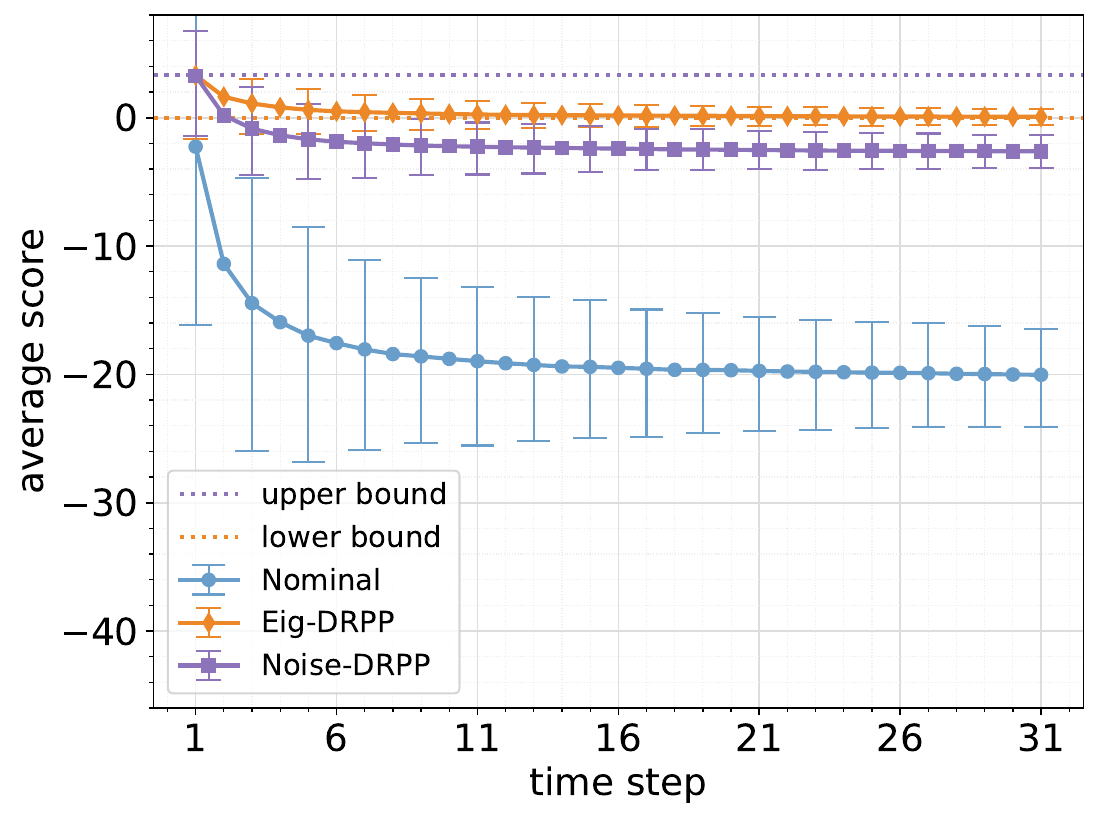}
    \caption{Prediction performance of different probabilistic predictors on the nonlinear manipulator system.}
    \label{fig:nonlinear}
\end{figure}

The nominal model adopts a lightweight numerical approximation that preserves the main physical structures while remaining efficient for repeated DRPP computations. For the inertia matrix, a diagonal-dominant rigid-body approximation is constructed by summing per-link COM contributions, i.e.,
$M_{i j}(q)=\sum_{k=\max (i, j)}^6 m_k l_{c, k}^2$.
For the Coriolis term, instead of computing full Christoffel symbols, we use a velocity-coupling surrogate, $C_{i j}(q, \dot{q})=k_c\left(\left|\dot{q}_i\right|+\left|\dot{q}_j\right|\right), k_c=10^{-3}$, which captures the dominant dependence on joint velocities.
The gravitational torque is computed from link COM heights, $g_i(q)=m_i g \sum_{k=1}^i l_{c, k} \cos (\sum_{j=1}^k q_j)$, where the gravitational acceleration $g=9.81$m/s$^2$.
Finally, the nominal discrete-time dynamics are obtained by RK4 integration with sampling time $\Delta t=0.1$s:
\begin{equation*}
    x_{k+1}=F\left(x_k, u_k\right)+w_k,
\end{equation*}
where $u_k=\tau(k\Delta t)$ is the zero-holder of the control input $\tau$ in the time interval $[k\Delta t, (k+1)\Delta t)$ and $w_k$ models the additive process noise in the state space. For the conditional concic moment-based ambiguity set, we let $\gamma_0(z) = \min\{0.3\|z\|_2, 5\}\Delta t^2$, $\gamma_1 = 0.5\Delta t^2$, $\gamma_2 = 3$, $\bar{\mu}_k = \mathbf{0}_{12\times 1}$, the nominal covariance for each pair $(q_i, \dot{q}_i)$ is $\begin{pmatrix}
        \Delta t^2 & 0.5 \Delta t^2\\
        0.5\Delta t^2 & 1.5\Delta t^2
    \end{pmatrix}$, and the full nominal covariance $\bar{\Sigma}_k$ is the blockwise combination of these $6$ sub-covariances. We simulate the nonlinear SDS under zero control by adversarially choosing one from the ambiguity set at each step against the predictor. The prediction performances of different probabilistic predictors are visualized and compared in Fig. \ref{fig:nonlinear}. Similar to the results in Fig. \ref{fig:performance}, both Noise-DRPP and Eig-DRPP outperform the nominal predictor significantly in terms of expected average score and score concentration. This experiment demonstrates the effectiveness and practicality of the proposed DRPP framework on nonlinear SDSs.

\section{Conclusion}
This paper presented a distributionally robust probabilistic prediction (DRPP) framework to optimize the worst-case prediction performance over a predefined ambiguity set of SDSs. To overcome the inherent intractability of optimizing over function spaces, Bellman equation is developed for the robust value function of DRPP, and optimality conditions are exploited to transform the original problem into Euclidean spaces.
While achieving the global optimal solution remains computationally prohibitive, we developed two suboptimal predictors: Noise-DRPP, obtained by relaxing the inner constraints, and Eig-DRPP, derived from relaxing the outer constraints. An explicit upper bound on the optimality gap was established for the proposed predictors.
Numerical simulations demonstrated the effectiveness and practicality of the proposed predictors across various SDSs. 

%\newpage  

\bibliographystyle{IEEEtran}
\bibliography{ref.bib}

\begin{appendices}
\section{Proof of Theorem \ref{thm:DPP}} \label{app:thm:DPP}
\begin{proof}
    For $k\in\{0, \ldots, T-1\}$ and any $z\in\mathcal{Z}$, we denote $W_k(z) = \sup_{\hat{p}_k \in \mathcal{F}} \inf_{\rho_k \in \mathcal{I}_k(z)} \int_{\mathcal{X}} \mathcal{L}(\hat{p}, x) + V^*_{k+1}(z^\prime) \mathrm{d} \rho_k(x)$. \\
    \underline{Step I: $V_k^*(z) \leq W_k(z)$.}
    According to the definitions \eqref{eq:obj},\eqref{eq:worst-obj}, and \eqref{eq:robust-value}, we have
    \begin{equation*}
        \begin{aligned}
            V_k^*(z) \!=& \!\!\sup_{\mathscr{F}_{k:T-1}} \inf_{\mathscr{P}_{k:T-1}} \!\!\mathbb{E}_{\mathscr{P}_{k:T-1}} \!\!\left[\sum_{t=k}^{T-1} \mathcal{L}(\mathscr{F}_t(\boldsymbol{z}_t), \boldsymbol{x}_{t+1}) \bigg| \boldsymbol{z}_k \!=\! z\right]\\
            \overset{(i)}{=}& \!\!\sup_{\mathscr{F}_{k:T-1}} \inf_{\mathscr{P}_{k:T-1}} \!\!\mathbb{E}_{\mathscr{P}_k}\bigg\{ \mathcal{L}(\mathscr{F}_k(z), \boldsymbol{x}_{k+1}) + \\
            & \qquad\quad\left.\mathbb{E}_{\mathscr{P}_{k+1:T-1}}\!\! \left[\sum_{t=k}^{T-1} \mathcal{L}(\mathscr{F}_t(\boldsymbol{z}_t), \boldsymbol{x}_{t+1})\right]\Bigg| \boldsymbol{z}_k \!=\! z\right\}\\
            \overset{(ii)}{=}& \!\!\sup_{\mathscr{F}_{k:T-1}} \inf_{\mathscr{P}_{k}} \mathbb{E}_{\mathscr{P}_k}\bigg\{ \mathcal{L}(\mathscr{F}_k(z), \boldsymbol{x}_{k+1}) + \\
            & \left.\inf_{\mathscr{P}_{k+1:T-1}}\!\!\mathbb{E}_{\mathscr{P}_{k+1:T-1}}\!\! \left[\sum_{t=k}^{T-1} \mathcal{L}(\mathscr{F}_t(\boldsymbol{z}_t), \boldsymbol{x}_{t+1})\right]\Bigg| \boldsymbol{z}_k \!=\! z\right\}\\
            =& \!\!\!\sup_{\mathscr{F}_{k:T-1}} \!\inf_{\mathscr{P}_{k}} \mathbb{E}_{\mathscr{P}_k}\!\!\left[ \mathcal{L}(\mathscr{F}_k(z), \boldsymbol{x}_{k+1}) \!+\! V_{k+1}^\mathscr{F}(\boldsymbol{z}_{k+1}) \bigg| \boldsymbol{z}_k \!=\! z\!\right]
        \end{aligned}
    \end{equation*}
    where $(i)$ holds because $\mathscr{P}_{k+1:T-1}$ do not affect the first term $\mathcal{L}(\mathscr{F}_k(z), \boldsymbol{x}_{k+1})$. According to Definition \ref{def:predictive_space} and \ref{def:ambiguity_set}, it follows that for all $\mathscr{F} \in \mathfrak{F}$, $\mathscr{P} \in \mathfrak{P}$, $z\in\mathcal{Z}$, and $t\in\{k, \ldots, T-1\}$, there is
    \begin{equation*}
        \begin{aligned}
            &\mathbb{E}_{\mathscr{P}_t}|\mathcal{L}(\mathscr{F}_t(\boldsymbol{z}_t), \boldsymbol{x}_{t+1}) \mid \boldsymbol{z}_t = z|\\
            \leq& \mathbb{E}_{\mathscr{P}_t} \left[C(1+\|\boldsymbol{x}_{t+1}\|_2^2)\mid \boldsymbol{z}_t=z\right] < \infty.
        \end{aligned}
    \end{equation*}
    Therefore, $(ii)$ interchanges $\inf_{\mathscr{P}_{k+1:T-1}}$ and $\mathbb{E}_{\mathscr{P}_k}$ according to the dominated convergence theorem.\\
    Notice that $V_{k+1}^\mathscr{F}(z) \leq V_{k+1}^*(z), \forall z\in\mathcal{Z}$, thus one has
    \begin{equation*}
        \begin{aligned}
            V_k^*(z) \leq& \sup_{\mathscr{F}_{k}} \inf_{\mathscr{P}_{k}} \mathbb{E}_{\mathscr{P}_k}\!\!\left[ \mathcal{L}(\mathscr{F}_k(z), \boldsymbol{x}_{k+1}) + V_{k+1}^*(\boldsymbol{z}_{k+1}) \bigg| \boldsymbol{z}_k \!=\! z\right]\\
            \overset{(i)}{=}& \sup_{\hat{p}_k \in \mathcal{F}} \inf_{\rho_k \in \mathcal{I}_k(z)} \int_{\mathcal{X}}\! \mathcal{L}(\hat{p}, x) + V^*_{k+1}(x, \pi_{k+1}(x)) \mathrm{d} \rho_k(x)\\
            =& W_k(z),
        \end{aligned}
    \end{equation*}
    where $(i)$ holds because: first, $\mathscr{F}_k$ affect the first term by $\mathcal{F}_k(z) \in \mathcal{F}$ and does not affect the second term $V_{k+1}^*(\boldsymbol{z}_{k+1})$; second, given $\boldsymbol{z}_k = z$, one has the conditional measure $P_k(\cdot \mid z) \in \mathcal{I}_k(z)$ according to Assumption \ref{ass:regularity}.\\
    \underline{Step II: $V_k^*(z) \geq W_k(z)$.}
    According to the definition $V_{k+1}^*(z) = \sup_{\mathscr{F}\in\mathfrak{F}} V_{k+1}^{\mathscr{F}}(z)$, it follows that for all $\epsilon > 0$, there exists a predictor $\mathscr{F}^\epsilon \in \mathfrak{F}$ such that $V_{k+1}^{\mathscr{F}^\epsilon}(z) \geq V_{k+1}^*(z) - \epsilon$ for all $z\in\mathcal{Z}$. Therefore,
    \begin{equation*}
        \begin{aligned}
            V_{k}^*(z) \!=&\!\!\! \sup_{\mathscr{F}_{k:T\!-\!1}}\! \inf_{\mathscr{P}_{k}} \mathbb{E}_{\mathscr{P}_k}\!\!\left[ \mathcal{L}(\mathscr{F}_k(z), \boldsymbol{x}_{k+1}) \!+\! V_{k+1}^\mathscr{F}(\boldsymbol{z}_{k+1})\bigg| \boldsymbol{z}_k \!=\! z\right]\\ 
            \geq& \sup_{\mathscr{F}_{k}} \inf_{\mathscr{P}_{k}} \mathbb{E}_{\mathscr{P}_k}\!\!\left[ \mathcal{L}(\mathscr{F}_k(z), \boldsymbol{x}_{k+1}) \!+\! V_{k+1}^{\mathscr{F}^\epsilon}(\boldsymbol{z}_{k+1})\bigg| \boldsymbol{z}_k \!=\! z\right]\\
            \geq& \sup_{\mathscr{F}_{k}} \inf_{\mathscr{P}_{k}} \mathbb{E}_{\mathscr{P}_k}\!\!\left[ \mathcal{L}(\mathscr{F}_k(z), \boldsymbol{x}_{k+1}) \!+\! V_{k+1}^{*}(\boldsymbol{z}_{k+1})\bigg| \boldsymbol{z}_k \!=\! z\right] \\
            &\qquad\qquad\qquad\qquad\qquad\qquad\qquad\qquad- \epsilon\\
            =& W_k(z) - \epsilon.
        \end{aligned}
    \end{equation*}
    Since $\epsilon > 0$ is arbitrary, it implies that $V_k^*(z) \geq W_k(z)$.\\    
    Combining the results from step $1$ and step $2$, we have $V_k^*(z) = W_k(z)$ for all $z\in\mathcal{Z}$.
    \end{proof}

\section{Duality of Problem \eqref{eq:p_0_cannolized}}\label{app:lem:lower-dual-problem}
\begin{proof}
The Lagrange function for \eqref{eq:p_0_cannolized} is defined as
\begin{equation*}
    \begin{aligned}
        & L(\rho_k, \nu_k,\mu_k, r, q, Q_1, Q_2, P_1, p_1, s_1, P_2, p_2, s_2)\\
        :=& \langle \rho_k, \log \hat{p}_k + V_{k+1}^*(x,\pi_{k+1}(x)) \rangle + r\left[\langle \rho_k, 1 \rangle \!-\! 1\right] \\
        &+q^\top\!\left[\langle \rho_k, x \rangle\!-\!\mu_k\!-\!\nu_k\right]\\
        &+ Q_1 \cdot \left[ \langle \rho_k, (x-\nu_k-\bar{\mu}_k)(x-\nu_k-\bar{\mu}_k)^\top \rangle - \gamma_{2} \bar{\Sigma}_k\right]\\
        &-Q_2 \cdot \left[ \langle \rho_k, (x-\nu_k-\bar{\mu}_k)(x-\nu_k-\bar{\mu}_k)^\top \rangle - \gamma_{3} \bar{\Sigma}_k\right]\\
        & - P_1\cdot\bar{\Sigma}_k - 2p_1^\top(\mu-\bar{\mu}_k) - s_1\gamma_1\\
        & - P_2\cdot I - 2p_2^\top(\nu_k-\bar{\nu}_k) - s_2\gamma_{0}(z).
    \end{aligned}
\end{equation*}
where the Lagrange multipliers satisfy that $Q_i \succeq 0$ and $\begin{bmatrix}
    P_i& p_i\\ p_i^\top & s_i
\end{bmatrix} \succeq 0$ for $i\in\{1,2\}$. Minimizing over $\rho_k \in \mathcal{M}_{+}(\mathcal{X})$, $\mu_k \in \mathbb{R}^{d_x}$, and $\nu_k\in\mathbb{R}^{d_x}$ we have
\begin{equation*}
    \begin{aligned}
        &\inf_{\rho_k, \mu_k,\nu_k} \!\!L(\rho_k, \nu_k,\mu_k, r, q, Q_1, Q_2, P_1, p_1, s_1, P_2, p_2, s_2)\\
        =&\inf_{\rho_k, \mu_k,\nu_k} \big\langle \rho_k, \log\hat{p}_k(x)+r+q^\top x+x^\top\!(Q_1\!-\!Q_2)x\\
        &\qquad\qquad+V_{k+1}^*(x,\pi_{k+1}(x)) \big\rangle - (q+2p_1)^\top\mu_k\\
        &- (q+2p_2)^\top\nu_k - (\gamma_{2}Q_1 - \gamma_{3}Q_2)\!\cdot\!\bar{\Sigma}_k \\
        &- (\nu_k+2\mu_k-\bar{\mu}_k)^\top(Q_1-Q_2)(\nu_k+\bar{\mu}_k)- r\\
        &- P_1\!\cdot\!\bar{\Sigma}_k+ 2p_1^\top\bar{\mu}_k- s_1\gamma_1 - P_2\!\cdot\! I+ 2p_2^\top\bar{\mu}_k- s_2\gamma_{0}(z) \\
        \overset{(i)}{=}&\inf_{\mu_k,\nu_k} -(q\!+\!2p_1)^\top\!\mu_k \!-\!(q\!+\!2p_2)^\top\nu_k \!-\! (\gamma_{2}Q_1 \!-\! \gamma_{3}Q_2)\!\cdot\!\bar{\Sigma}_k  \\
        &- (\nu_k+2\mu_k-\bar{\mu}_k)^\top(Q_1-Q_2)(\nu_k+\bar{\mu}_k)- r\\
        &- P_1\!\cdot\!\bar{\Sigma}_k+ 2p_1^\top\bar{\mu}_k- s_1\gamma_1 - P_2\!\cdot\! I+ 2p_2^\top\bar{\nu}_k- s_2\gamma_{0}(z) \\
        \overset{(ii)}{=}&\inf_{\nu_k}\!-r\!+\! \bar{\mu}_k^\top\!(Q_1\!-\!Q_2)\bar{\mu}_k\!-\! (\gamma_{2}Q_1 \!-\! \gamma_{3}Q_2\!+\!P_1)\!\cdot\!\bar{\Sigma}_k\! - P_2\cdot I \\
        &+\! 2p_1^\top\!\bar{\mu}_k\!-\! s_1\gamma_1+\! 2p_2^\top\!\bar{\mu}_k\!-\! s_2\gamma_{0}(z) - (q+2p_2)^\top\nu_k\\
        &- \nu_k^\top(Q_1-Q_2)\nu_k,
    \end{aligned}
\end{equation*}
where $q+2p_1+2(Q_1-Q_2)(\nu_k+\bar{\mu}_k) = 0$.
Equation $(i)$ holds when $\log\hat{p}_k(x)+r+q^\top x+x^\top(Q_1-Q_2)x+V_{k+1}^*(x,\pi_{k+1}(x)) \geq 0 \;\forall x\in\mathbb{R}^{d_x}$, otherwise there is $\inf_{\rho_k \in \mathcal{M}_{+}(\mathcal{X})}L$ equals $-\infty$. Equation $(ii)$ holds because
\begin{equation*}
    \begin{aligned}
        &\inf_{\mu_k,\nu_k}- (q+2p_1)^\top\mu_k -(q+2p_2)^\top\nu_k-\bar{\mu}_k^\top(Q_1-Q_2)\bar{\mu}_k\\
        &-(\nu_k+2\mu_k-\bar{\mu}_k)^\top(Q_1-Q_2)(\nu_k+\bar{\mu}_k)\\
        =& \inf_{\nu_k}\inf_{\mu_k}-[q+2p_1+2(Q_1-Q_2)(\nu_k+\bar{\mu}_k)]^\top \mu_k  \\
        &\qquad - (q+2p_2)^\top\nu_k- \nu_k^\top(Q_1-Q_2)\nu_k \\
        =& \inf_{\nu_k} - (q+2p_2)^\top\nu_k - \nu_k^\top(Q_1-Q_2)\nu_k\\
        &\text{ s.t. } q+2p_1+2(Q_1-Q_2)(\nu_k+\bar{\mu}_k) = 0.\\
        % =& \left[\!\frac{1}{2}q\!+\!2p_2\!-\!p_1\!-\!(Q_1\!-\!Q_2)\bar{\mu}_k\right]^{\!\!\top\!}\!\! \left[\frac{1}{2}(Q_1\!-\!Q_2)^{-1}\!(q\!+\!2p_1)\!+\!\bar{\mu}_k\!\right]
    \end{aligned}
\end{equation*}
Combining those conditions and the previous multipliers' constraints, we have the dual form completed.
\end{proof}

\section{Proof of Theorem \ref{thm:opt-condition}}\label{app:thm:opt-condition}
\begin{proof}
\underline{Step I.} We leverage the separable structure of $\mathbf{D}_1$ to reformulate the problem. For ease of notation, we use the $\kappa$ to summarize all the Lagrange multipliers except for $r$, i.e.,
\[
    \kappa = (q,Q_1,Q_2,P_1,p_1,s_1,P_2,p_2,s_2).
\]
First, the objective function can be separated as $G(r,\kappa) = -r + g(\kappa)$ where $g(\kappa) = \bar{\mu}_k^\top\!(Q_1\!-\!Q_2)\bar{\mu}_k- (\gamma_{2}Q_1 - \gamma_{3}Q_2+P_1)\cdot\bar{\Sigma}_k- P_2\cdot I + 2p_1^\top\bar{\mu}_k- s_1\gamma_1+2p_2^\top\bar\nu_k - s_2\gamma_{0}(z) - (q+2p_2)^\top\nu_k^* - \nu_k^{*,\top}(Q_1-Q_2)\nu_k^*$. Second, only one of the constraints concerns $r$. Therefore, given any group of feasible $\kappa$ satisfying the other constraints in $\mathbf{D}_1$, the optimizing over $\hat{p}_k$ and $r$ is equivalent to solving the following problem:
\begin{equation}\label{eq:separate-sup}
    \begin{aligned}
        &\sup_{\hat{p}_k\in\mathcal{F},r} -r+g(\kappa)\\
        &\text{s.t.} \left\{
        \begin{aligned}
            &x^\top(Q_1-Q_2)x + x^\top q + r + \log \hat{p}_k(x)\\
            &+ V_{k+1}^*(x,\pi_{k+1}(x))\geq 0 \; \forall x \in \mathcal{X}
        \end{aligned}
        \right.
    \end{aligned}
\end{equation}

\underline{Step II.} We exploit the last constraint in \eqref{eq:separate-sup} to analyze the lower bound of $r$. Notice that
\begin{equation}\label{eq:constraint}
    \begin{aligned}
        1 &= \int_{\mathbb{R}^{d_x}} \hat{p}_k(x)\mathrm{d}x \\
        &\geq \int_{\mathbb{R}^{d_x}} \exp\big\{-x^\top\!(Q_1-Q_2)x - x^\top q - r\\
        &\qquad\qquad\qquad -V_{k+1}^*(x,\pi_{k+1}(x))\big\} \mathrm{d}x,
    \end{aligned}
\end{equation}
one immediately gets that $r$ is lower bounded by $r^{*} := \log[\int_{\mathbb{R}^{d_x}} \!\!\exp\{-x^\top\!(Q_1-Q_2)x - x^\top q - V_{k+1}^*(x,\pi_{k+1}(x))\} \mathrm{d}x]$. 

Notice that the equality of \eqref{eq:constraint} holds only when there is $\hat{p}_k(x) = \exp\{-x^\top\!(Q_1-Q_2)x - x^\top q - r^*\}$ for $x \in \mathbb{R}^{d_x}$ almost everywhere. Therefore, given any feasible $\kappa$, the objective function $-r+g(\kappa)$ in problem \eqref{eq:separate-sup} can only attain its maximum when $\hat{p}_k$ belongs to the following exponential family almost everywhere on $\mathbb{R}^{d_x}$, 
\begin{equation*}
    \hat{p}_k(x) \overset{a.s.}{\propto} \exp\{-x^\top \theta_1 x - x^\top \theta_2 - V_{k+1}^*(x, \pi_{k+1}(x))\},
\end{equation*}
where $\theta_1 \in S^{d_x}$ and $\theta_2 \in \mathbb{R}^{d_x}$.

\underline{Step III.}
Finally, when $k=T-1$, there is $V_{T}^*(x,\pi_{T}(x)) = 0$. Because $Q_1$ and $Q_2$ are both semi-definite, we have $Q_1-Q_2$ is symmetric and there exists an orthogonal matrix $U\in\mathbb{R}^{d_x\times d_x}$ and an diagonal matrix $D = \operatorname{diag}(\lambda_1,\ldots, \lambda_{d_x})$ such that $Q_1-Q_2 = UDU^\top$. Then we have
\begin{equation*}
    \begin{aligned}
        &\int_{\mathbb{R}^{d_x}} \exp\{-x^\top(Q_1-Q_2)x - x^\top q - r\} \mathrm{d}x\\
        =& \int_{\mathbb{R}^{d_x}} \exp\{-x^\top UDU^\top x - x^\top q - r\} \mathrm{d}x\\
        % =& \int_{\mathbb{R}^{d_x}} \exp\{-(U^\top x)^\top DU^\top x \!-\! (U^\top x)^\top U^\top q \!-\! r\} |U|\mathrm{d}(U^\top x)\\
        =& \int_{\mathbb{R}^{d_x}} \exp\{-\tilde{x}^\top D\tilde{x} - \tilde{x}^\top \tilde{q} - r\} \operatorname{det}(U)\mathrm{d}\tilde{x}\\
        =&e^{-r}\operatorname{det}(U) \prod_{i=1}^{d_x} \int_{\mathbb{R}} \exp\{-\tilde{x}_{(i)}^2\lambda_i - \tilde{x}_{(i)}\tilde{q}_{(i)}\} \mathrm{d}\tilde{x}_{(i)},
    \end{aligned}
\end{equation*}
where the second equation follows by letting $\tilde{x} = U^\top x$ and $\tilde{q} =U^{\top} q$. If there is an index $j\in\{1,\ldots,d_x\}$ such that $\lambda_j < 0$, then $\int_{\mathbb{R}} \exp\{-\tilde{x}^2_{(i)}\lambda_j - \tilde{x}_{(j)}\tilde{q}_{(j)}\}\mathrm{d}\tilde{x}_{(j)} = \infty$, which contradicts \eqref{eq:constraint}. 
Therefore, $Q_1-Q_2$ is semi-definite positive, and a necessary optimal condition for the inner optimization is that $\hat{p}_{T-1}$ is subject to Gaussian almost everywhere. 
\end{proof}

\section{Derivation of \eqref{eq:p_0_cannolized-dual-join-parametrize}}\label{app:use_optimality_condition}
\begin{proof}
Facilitated by Theorem \ref{thm:opt-condition}, when $k=T-1$ we can first parameterize $\hat{p}_k$ by $\hat{\mu}_k\in\mathbb{R}^{d_x}$ and $\hat{\Sigma}_k\in S_{+}^{d_x}$ such that
\begin{equation*}
    \begin{aligned}
        \hat{p}_k(x) \!=\! 
        \frac{1}{\sqrt{(2\pi)^{d_x}\operatorname{det}(\hat\Sigma_k)}}\exp\!\left[\!-\frac{1}{2}(x\!-\!\hat\mu_k)^\top\! {\hat\Sigma_k}^{-1}(x\!-\!{\hat\mu_k})\right].
    \end{aligned}
\end{equation*}
Substituting it into the constraints of $\mathbf{D}_1$, there is
\begin{equation*}
    \begin{aligned}
        &x^\top(Q_1-Q_2)x + x^\top q + r -\frac{1}{2}\left[d_x\log(2\pi)+\log\operatorname{det}(\hat\Sigma_k)\right] \\
        &-\frac{1}{2}(x-\hat\mu_k)^\top {\hat\Sigma_k}^{-1}(x-\hat\mu_k) = 0 \quad \forall x \in \mathbb{R}^{d_x}.
    \end{aligned}
\end{equation*}
Then we have $Q_1 - Q_2 = \frac{1}{2} \hat{\Sigma}_k^{-1}$, $q = -\hat{\Sigma}^{-1}_k\hat{\mu}_k$, and $r = \frac{1}{2}\left[d_x\log(2\pi)+\log\operatorname{det}(\hat\Sigma_k)+\hat{\mu}_k^\top\hat{\Sigma}_k^{-1}\hat{\mu}_k\right]$. Substituting these equations and the parameterized $\hat{p}_k$ into the objective of problem $\mathbf{D}_1$, we get 
\begin{equation*}
    \begin{aligned}
        &G(r,q,Q_1,Q_2,P_1,p_1,s_1,P_2,p_2,s_2)\\
        % =& - r \!+\! \bar{\mu}_k^\top\!(Q_1\!-\!Q_2)\bar{\mu}_k- (\gamma_{2}Q_1 - \gamma_{3}Q_2+P_1)\cdot\bar{\Sigma}_k\\
        % &- P_2\cdot I + 2p_1^\top\bar{\mu}_k- s_1\gamma_1+2p_2^\top\bar\nu_k - s_2\gamma_{0}(z)+ \\
        % &[\frac{1}{2}q+2p_2-p_1-(Q_1-Q_2)\bar{\mu}_k]^{\top}\\
        % &[\frac{1}{2}(Q_1-Q_2)^{-1}(q+2p_1)+\bar{\mu}_k]\\
        =& -\frac{1}{2}\left[d_x\log(2\pi)+\log\operatorname{det}(\hat\Sigma_k)\right]\\
        &- P_2\cdot I - s_1\gamma_1+2p_2^\top(\bar\nu_k+\bar{\mu}_k-\hat{\mu}_k) - s_2\gamma_{0}(z)\\
        & + 2(2p_2-p_1)^\top\hat{\Sigma}_k p_1- (\gamma_{2}Q_1\!-\!\gamma_{3}Q_2\!+\!P_1)\!\cdot\!\bar{\Sigma}_k\\
        =& \!-\!\frac{1}{2}d_x\log(2\pi) \!+\! \frac{1}{2}\log\operatorname{det}(\hat{\Sigma}_k^{-1}) \!-\! (\gamma_{2}Q_1\!-\!\gamma_{3}Q_2)\!\cdot\!\bar{\Sigma}_k\\
        &-P_1\!\cdot\!\bar{\Sigma}_k - P_2\!\cdot\! I - s_1\gamma_1 - s_2\gamma_{0}(z)\\
        &- 2p_2^\top(\hat{\mu}_k-\bar{\mu}_k-\bar{\nu}_k)+ 2p_1^\top\!\hat{\Sigma}_k (2p_2\!-\!p_1).
    \end{aligned}
\end{equation*}
Substituting the above equations into the problem $\mathbf{D}_1$, we have the problem \eqref{eq:p_0_cannolized-dual-join-parametrize} derived.
\end{proof}

\section{Proof of Lemma \ref{lem:primal-solution-1}}\label{app:lem:primal-solution-1}
\begin{proof}
    Since the objective function is linear with respect to $\Sigma_k$ and the second constraint provides an upper bound for it, one immediately has the worst-case covariance should satisfy that
    \begin{equation*}
        \Sigma_k^{*} = \gamma_{2}\bar{\Sigma}_k - (\mu_k-\bar{\mu}_k)(\mu_k-\bar{\mu}_k)^\top.
    \end{equation*}
    
    Let $a = \mu_k-\bar{\mu}_k, b = \nu_k-\bar{\nu}_k, c = \bar{\mu}_k+\bar{\nu}_k-\hat{\mu}_k$, problem $\mathbf{P}_2$ can be further equivalent to
    \begin{equation}\label{eq:one-step-DRPP-general-parameterized-2}
        \begin{aligned}
            \inf_{\hat{\Sigma}_k^{-1}, c}\sup_{a,b} &\;-\log\operatorname{det}(\hat{\Sigma}_k^{-1})+\gamma_{2}\bar{\Sigma}_k\!\cdot\!\hat{\Sigma}_k^{-1} \\
            &\; +\|a+b+c\|_{\hat{\Sigma}_k^{-1}}^2 - \|a\|_{\hat{\Sigma}_k^{-1}}^2\\
            \text{s.t.}&\;            \|a\|_{\bar{\Sigma}_k^{-1}}^2\leq \gamma_1, \|b\|_2^2 \leq \gamma_{0}(z).
        \end{aligned}
    \end{equation}
    Suppose $(a^{*}, b^{*}) \in \arg\max_{a,b}\left\{ \|a+b\|_{\hat{\Sigma}^{-1}}^2 - \|a\|_{\hat{\Sigma}^{-1}}^2\right\}$, we have $(-a^{*},-b^{*})$  can also attain the maximum. Then we have
    \begin{equation*}
        \begin{aligned}
            &\max_{a,b}\|a+b+c\|^2_{\hat{\Sigma}_k^{-1}} - \|a\|^2_{\hat{\Sigma}_k^{-1}}\\
            \geq& \max\left\{\|a^{*}+b^{*}+c\|^2_{\hat{\Sigma}_k^{-1}} - \|a^{*}\|^2_{\hat{\Sigma}_k^{-1}},\right. \\
            &\left.\qquad\|-a^{*}-b^{*}+c\|^2_{\hat{\Sigma}_k^{-1}} - \|-a^{*}\|^2_{\hat{\Sigma}_k^{-1}}\right\}\\
            =&\|a^{*}+b^{*}\|^2_{\hat{\Sigma}_k^{-1}} - \|a^{*}\|^2_{\hat{\Sigma}_k^{-1}} + \|c\|_{\hat{\Sigma}_k}^{-1} \\
            &+ 2\max\left\{\langle c,a^{*}+b^{*} \rangle_{\hat{\Sigma}_k^{-1}}, -\langle c,a^{*}+b^{*} \rangle_{\hat{\Sigma}_k^{-1}}\right\}\\
            \geq& \|a^{*}+b^{*}\|^2_{\hat{\Sigma}_k^{-1}} - \|a^{*}\|^2_{\hat{\Sigma}_k^{-1}}
        \end{aligned}
    \end{equation*}
    where the last inequality is strict when $c\neq 0$. Therefore, we have $c^{*} = \arg\min_{c} \max_{a,b}\|a+b+c\|^2_{\hat{\Sigma}_k^{-1}} - \|a\|^2_{\hat{\Sigma}_k^{-1}} = 0$, which means that $\hat{\mu}^{*}_k = \bar{\mu}_k + \bar{\nu}_k$.
    
% \section{Proof of Lemma \ref{lem:QCLP}}\label{app:lem:QCLP}
% Using the KKT condition, there are
% $$\left\{\begin{aligned}
% &\partial_{a} \left[\|a+b\|_{\hat{\Sigma}_k^{-1}}^2 - \|a\|_{\hat{\Sigma}_k^{-1}}^2 - s_1(\|a\|_{\bar{\Sigma}_k^{-1}}^2 - \gamma_1)\right] = 0\\
% & s_1(\|a\|_{\bar{\Sigma}_k^{-1}}^2 - \gamma_1) = 0, s_1\geq 0\\
% \end{aligned}\right.$$
% The first condition reveals that $s_1 \neq 0$ and $a = s_1^{-1}\bar{\Sigma}_k\hat{\Sigma}_k^{-1}b$, thus the second condition leads to $\|a\|_{\bar{\Sigma}_k^{-1}}^2 = \gamma_1$. Finally, we have
% $$a^{*} = \frac{\sqrt{\gamma_1}\bar\Sigma_k\hat{\Sigma}_k^{-1}b}{\|\bar\Sigma_k\hat{\Sigma}_k^{-1}b\|_{\bar{\Sigma}_k^{-1}}}.$$
\end{proof}

\section{Duality of the moment problem \eqref{eq:noise-DRPP}}\label{app:lem:upper-dual-problem}
\begin{proof}
The Lagrange function for \eqref{eq:noise-DRPP} is defined as
\begin{equation*}
    \begin{aligned}
        & L(\tilde{\rho}_k, \mu_k, r, q, Q_1, Q_2, P, p, s)\\
        :=& \langle \tilde{\rho}_k, \log \hat{p}_k(w\!+\!\bar\nu_k)\rangle \!+\! r\left[\langle \tilde{\rho}_k, 1 \rangle \!-\! 1\right] \!+\! q^\top\!\left[\langle \tilde{\rho}_k, w \rangle\!-\!\mu_k\right]\\
        &+ Q_1 \cdot \left[ \langle \tilde{\rho}_k, (w-\bar{\mu}_k)(w-\bar{\mu}_k)^\top \rangle - \gamma_{2} \bar{\Sigma}_k\right]\\
        &-Q_2 \cdot \left[ \langle \tilde{\rho}_k, (w-\bar{\mu}_k)(w-\bar{\mu}_k)^\top \rangle - \gamma_{3} \bar{\Sigma}_k\right]\\
        & - P\cdot\bar{\Sigma}_k - 2p^\top(\mu_k-\bar{\mu}_k) - s\gamma_1,
    \end{aligned}
\end{equation*}
where the Lagrange multipliers satisfy $Q_1 \succeq 0, Q_2 \succeq 0$ and $\begin{bmatrix}
    P& p\\ p^\top & s
\end{bmatrix} \succeq 0$. Minimizing over $\tilde{\rho}_k \in \mathcal{M}_{+}(\mathcal{X})$ and $\mu_k \in \mathbb{R}^{d_x}$, one has
\begin{equation*}
    \begin{aligned}
        &\inf_{\tilde{\rho}_k \in \mathcal{M}_{+}(\mathcal{X}), \mu_k} L(\tilde{\rho}_k, \mu_k, r, q, Q_1, Q_2, P, p, s)\\
        =&\inf_{\tilde{\rho}_k \in \mathcal{M}_{+}(\mathcal{X}), \mu_k} \!\!\!\!\langle \tilde{\rho}_k, \log\hat{p}_k(w\!+\!\bar\nu_k)\!+\!r\!+q\!^\top\!w\!+\!w^\top\!(Q_1\!-\!Q_2)w \rangle\\
        &-\left[q+2(Q_1-Q_2)\bar{\mu}_k+2p\right]^\top\mu_k - r+ \bar{u}_k^\top(Q_1-Q_2)\bar{u}_k\\
        &- (\gamma_{2}Q_1 - \gamma_{3}Q_2+P)\!\cdot\!\bar{\Sigma}_k+ 2p^\top\bar{\mu}_k- s\gamma_1 \\
        =&- r \!+\! \bar{u}_k^\top\!(Q_1\!-\!Q_2)\bar{u}_k\!-\! (\gamma_{2}Q_1 \!-\! \gamma_{3}Q_2\!+\!P)\!\cdot\!\bar{\Sigma}_k\!+\! 2p^\top\!\bar{\mu}_k\\
        &-\! s\gamma_1,
    \end{aligned}
\end{equation*}
where the second equality holds when $\log\hat{p}_k(w+\bar\nu_k)+r+q^\top w+w^\top(Q_1-Q_2)w \geq 0 \;\forall w\in\mathbb{R}^{d_x}$ and $q+2(Q_1-Q_2)\bar{\mu}_k+2p=0$ hold simultaneously, otherwise there is $\inf_{\tilde{\rho}_k \in \mathcal{M}_{+}(\mathcal{X}), \mu_k}L = -\infty$. Combining these two conditions and the previous multipliers' constraints, we have attained the dual problem $\mathbf{D}_2$.
\end{proof}

\section{Proof of Theorem \ref{thm:one-step-DRPP-mp}}\label{app:thm:one-step-DRPP-mp}
\begin{proof}
The proof is organized in three parts: first, we use Theorem \ref{thm:opt-condition} to parameterize $\hat{p}_k$ by Gaussian distribution $\mathcal{N}(\hat{\mu}_k, \hat{\Sigma}_k)$ and reformulate the dual problem $\mathbf{D}_2$ to a finite-dimensional optimization; second, we solve the optimal predictive pdf $\hat{p}^{*}_k\sim\mathcal{N}(\hat{\mu}^{*}_k, \hat{\Sigma}^{*}_k)$; third, we substitute $\hat{p}^{*}_k$ into the original problem $\mathbf{P}_1$ to obtain the worst-case conditional measure $\rho^{*}_k$ by solving the inner minimization.

\underline{Step I.} Using the necessary optimality condition for $\hat{p}_k$, we can parameterize $\hat{p}_k$ by $\hat\mu_k\in\mathbb{R}^{d_x}$ and $\hat\Sigma_k\in S_{+}^{d_x}$ such that
\begin{equation*}
    \begin{aligned}
        \hat{p}_k(w\!+\!\bar{\nu}_k) \!=\! 
        \frac{\exp\!\left[\!-\frac{1}{2}(w\!-\!\hat\mu_k)^\top\! {\hat\Sigma_k}^{-1}(w\!-\!{\hat\mu_k})\right]}{\sqrt{(2\pi)^{d_x}\operatorname{det}(\hat\Sigma_k)}}.
    \end{aligned}
\end{equation*}
Moreover, from the constraint $w^\top\!(Q_1-Q_2)w + w^\top q + r + \log \hat p_k(w+\bar\nu_k) \geq 0 \; \forall w \in \mathbb{R}^{d_x}$, we get
\begin{equation*}
    \begin{aligned}
        &w^\top(Q_1\!-\!Q_2)w + w^\top\! q + r -\frac{1}{2}\left[d_x\log(2\pi)\!+\!\log\operatorname{det}(\hat\Sigma_k)\right] \\
        &-\frac{1}{2}(w-\hat\mu_k)^\top {\hat\Sigma_k}^{-1}(w-\hat\mu_k) = 0 \quad \forall w \in \mathbb{R}^{d_x},
    \end{aligned}
\end{equation*}
which indicates that $Q_1 - Q_2 = \frac{1}{2} \hat{\Sigma}_k^{-1}$, $q = -\hat{\Sigma}_k^{-1}\hat{\mu}_k$, and $r = \frac{1}{2}\left[d_x\log(2\pi)+\log\operatorname{det}(\hat\Sigma_k)+\hat{\mu}_k^\top\hat{\Sigma}_k^{-1}\hat{\mu}_k\right]$. Next, from the constraint $q + 2(Q_1-Q_2)\bar{\mu}_k + 2p = 0$ one can get $p = \frac{1}{2}\hat{\Sigma}_k^{-1}(\hat{\mu}_k-\bar{\mu}_k)$. 
Substituting the parameterized $\hat{p}_k$ and the above equations into problem $\mathbf{D}_2$, there is 
\begin{equation*}
    \begin{aligned}
        &G(r,q,Q_1,Q_2,P,p,s) = -\frac{1}{2}\left[d_x\log(2\pi)+\log\operatorname{det}(\hat\Sigma_k)\right]\\
        &\quad\;-\frac{1}{2}\hat{\mu}_k^\top\! {\hat\Sigma_k}^{-1}\hat{\mu}_k+\bar{\mu}_k^\top(Q_1-Q_2)\bar{\mu}_k + \bar{\mu}_k^\top\hat{\Sigma}_k^{-1}(\hat{\mu}_k-\bar{\mu}_k)\\
        &\quad\;- (\gamma_{2}Q_1-\gamma_{3}Q_2+P)\cdot\bar{\Sigma}_k - s\gamma_1\\
    &= -\frac{1}{2}\left[d_x\log(2\pi) - \log\operatorname{det}(\hat{\Sigma}_k^{-1}) \!+\!  \|\hat\mu_k\!-\!\bar{\mu}_k\|_{\hat{\Sigma}_k^{-1}}^2\right]\\
        &\quad\;- (\gamma_{2}Q_1-\gamma_{3}Q_2+P)\cdot\bar{\Sigma}_k - s\gamma_1.
    \end{aligned}
\end{equation*}
As a result, $\mathbf{D}_2$ can be transformed into a finite-dimensional optimization problem as follows:
\begin{equation}\label{eq:noise-DRPP-dual-join-parametrize}
    \begin{aligned}
        &\sup_{\hat{\mu}_k,\hat{\Sigma}_k,Q_1,Q_2,P,p,s}\!\! -\frac{1}{2}d_x\log(2\pi) \!+\! \frac{1}{2}\log\operatorname{det}(\hat{\Sigma}_k^{-1})\!-\! s\gamma_1\\
        &-\frac{1}{2}\|\hat\mu_k\!-\!\bar{\mu}_k\|_{\hat{\Sigma}_k^{-1}}^2\!-\! (\gamma_{2}Q_1 \!-\! \gamma_{3}Q_2\!+\!P)\!\cdot\!\bar{\Sigma}_k \\          
        &\text{s.t.} \left\{
        \begin{aligned}
            &Q_1 \succeq 0, Q_2 \succeq 0, Q_1 - Q_2 = \frac{1}{2}\hat{\Sigma}_k^{-1}, \hat{\Sigma}_k \succ 0\\
            &\left[\begin{array}{cc}
                P & p\\
                p^\top & s
            \end{array}\right]\succeq 0, p = \frac{1}{2}\hat{\Sigma}_k^{-1}(\hat{\mu}_k-\bar{\mu}_k).
        \end{aligned}
        \right.
    \end{aligned}
\end{equation}

\underline{Step II.} Notice that the constraint $\hat{\Sigma}_k\succ 0$ and the relationship $p = \frac{1}{2}\hat{\Sigma}_k^{-1}(\hat{\mu}_k-\bar{\mu}_k)$ reveals that $p$ and $\hat{\mu}_k$ are one-on-one, one can optimize on either of them. In this proof, we choose to maximize over $p$ and substitute $\hat{\mu}_k$ in the objective by $p$. Next, notice that the constraints for $Q_1, Q_2$ only depend on $\hat{\Sigma}_k$, and it can be separated from the constraints for $P, p, s$, we can first maximize over $Q_1$ and $Q_2$, then over $P, p, s$, finally over $\hat{\Sigma}$. In other words, the objective of problem \eqref{eq:noise-DRPP-dual-join-parametrize} can be reformulated as
\begin{equation}\label{eq:noise-DRPP-dual-join-parametrize-re}
    \begin{aligned}
        &\sup_{\hat{\Sigma}_k}\sup_{P,p,s}\sup_{Q_1,Q_2}\!\! -\frac{1}{2}d_x\log(2\pi) \!+\! \frac{1}{2}\log\operatorname{det}(\hat{\Sigma}^{-1})\\
        &-2p^\top\!\hat{\Sigma}p\!-\! (\gamma_{2}Q_1 \!-\! \gamma_{3}Q_2\!+\!P)\!\cdot\!\bar{\Sigma}_k-\! s\gamma_1.
    \end{aligned}
\end{equation}

The solution of the first maximization is $Q_1^{*} = \frac{1}{2}\hat{\Sigma}_k^{-1}, Q_2^{*} = 0$ because
\begin{equation*}
\begin{aligned}
    Q_1^{*} =& \arg\max_{Q_1\succeq\frac{1}{2}\hat{\Sigma}_k^{-1}} -((\gamma_{2}-\gamma_{3})Q_1 + \frac{\gamma_{3}}{2}\hat{\Sigma}_k^{-1}) \cdot \bar{\Sigma}_k\\
    =& \arg\max_{Q_1\succeq\frac{1}{2}\hat{\Sigma}_k^{-1}} -(\gamma_{2}-\gamma_{3})Q_1 \cdot \bar{\Sigma}_k= \frac{1}{2}\hat{\Sigma}_k^{-1},
\end{aligned}
\end{equation*}
and $Q_2^{*} = Q_1^{*} - \frac{1}{2}\hat{\Sigma}_k^{-1} = 0$.

The solution of the second maximization is $P^{*} = 0, p^{*} = 0, s^{*} = 0$. Notice that $\begin{bmatrix}
    P & p\\ p^\top& s
\end{bmatrix}\succeq 0$ is equivalent to $s \geq p^\top P^{-1} p$, therefore $s\geq 0$ and $\arg\max_{s} - s\gamma_1 = 0$. Again, $\arg\max_{p} -2p^\top\hat{\Sigma}_kp = 0$ and $\arg\max_{P} P\cdot\bar{\Sigma}_k = 0$. Additionally, $p^{*} = 0$ indicates that $\hat{\mu}^{*}_k = 2\hat{\Sigma}_kp + \bar{\mu}_k = \bar{\mu}_k$.

The solution of the third maximization is $\hat{\Sigma}^{*}_k = \gamma_{2}\bar{\Sigma}_k$. This is because $\hat{\Sigma}^{*}_k = [\arg\max_{\hat{\Sigma}_k^{-1}\succ 0} F(\hat{\Sigma}_k^{-1})]^{-1}$ where $F(X) = \log\operatorname{det}(X) - X \cdot \gamma_{2}\bar{\Sigma}_k$. Notice that $\frac{\mathrm{d}}{\mathrm{d}X}F(X) = X^{-1} - \gamma_{2}\bar{\Sigma}_k$ and $\frac{\mathrm{d}}{\mathrm{d}X}F$ monotonously decrease on $X$. Therefore, $F$ is concave, which can be maximized at $X^{*}$ where $\frac{\mathrm{d}}{\mathrm{d}X}F(X^{*}) = 0$, i.e., $X^{*} = (\gamma_{2}\bar{\Sigma}_k)^{-1}$. Finally, we have $\hat{\Sigma}^{*}_k$ equals the inverse of $X^{*}$, which is $\gamma_{2}\bar{\Sigma}_k$.
Therefore, the optimal predictive pdf is $\hat{p}^{*}_k \sim \mathcal{N}\left(\bar{\nu}_k + \bar{\mu}_k, \gamma_{2}\bar{\Sigma}_k\right)$.

\underline{Step III.} Substituting the optimal $\hat{p}^{*}_k$ into problem $\mathbf{P}_1$, we have the inner minimization problem as 
\begin{equation}\label{eq:one-step-DRPP-worst}
    \begin{aligned}
        &\inf_{\tilde{\rho}_k,\mu_k} \; \mathbb{E}_{w \sim \tilde{\rho}_k} -\frac{1}{2}\left[d_x\log(2\pi)+\log\operatorname{det}(\gamma_{2}\bar{\Sigma}_k)\right.\\
        &\qquad\qquad\qquad\left.+(w-\bar{\mu}_k)^\top (\gamma_{2}\bar{\Sigma}_k)^{-1}(w-\bar{\mu}_k)\right]\\
        &\text{s.t.} \left\{
        \begin{aligned}
            & \mathbb{E}_{w\sim \tilde{\rho}_k}[1]=1, \mathbb{E}_{w\sim \tilde{\rho}_k}[w]=\mu_k  \\
            & \gamma_{3}\bar{\Sigma}_k \preceq \mathbb{E}_{w\sim \tilde{\rho}_k}\!\left[\left(w\!-\!\bar{\mu}_k\right)\left(w\!-\!\bar{\mu}_k\right)^{\mathrm{T}}\right] \preceq \gamma_{2} \bar{\Sigma}_k \\
            & \left[\begin{array}{cc}
            \bar{\Sigma}_k & \left(\mu_k-\bar{\mu}_k\right) \\
            \left(\mu_k-\bar{\mu}_k\right)^{\top} & \gamma_1
            \end{array}\right] \succeq 0.
        \end{aligned}
        \right.
    \end{aligned}
\end{equation}
The objective is minimized when $\mathbb{E}\!\left[\left(\boldsymbol{w}_k\!-\!\bar{\mu}_k\right)\left(\boldsymbol{w}_k\!-\!\bar{\mu}_k\right)^{\mathrm{T}}\right]$ equals $\gamma_{2} \bar{\Sigma}_k$, and the proof is completed.
% which means the solution for $p_{\boldsymbol{w}_k}^{*}$ is the set
% $$
% \left\{p_{\boldsymbol{w}_k} \mid  \mathbb{E}_{w\sim p_{\boldsymbol{w}_k}}\!\left[\left(w\!-\!\bar{\mu}_k\right)\left(w\!-\!\bar{\mu}_k\right)^{\mathrm{T}}\right] = \gamma_{2} \bar{\Sigma}_k \right\},
% $$
% and the objective function at $(\hat{p}^{*}_k, \rho^{*}_k)$ is
% \begin{equation*}
%     -\frac{1}{2}\left[d_x\log(2\pi)+d_x+\log\operatorname{det}(\gamma_{2}\bar{\Sigma}_k)\right].
% \end{equation*}
\end{proof}

% \section{Proof of Lemma \ref{lem:convex-max-bstar}}\label{app:lem:convex-max-bstar}
% There are necessary optimality conditions for \eqref{eq:convex-max}:
% \begin{equation*}
%     \left\{\begin{aligned}
%         & \left(\frac{\alpha}{\|b\|_A} A + 2B -2sI\right)b = 0\\
%         & \|b\|_2^2 = \gamma_{0}(z_k),
%     \end{aligned}\right.
% \end{equation*}
% where the first one comes from the KKT condition and the second one holds because the objective function increases monotonously with $\|b\|_2$.
% It follows that the optimal $b$ should be an eigenvector of $\frac{\alpha}{\|b\|_A} A + 2B$ with its $\ell^2$-norm being $\sqrt{\gamma_{0}(z_k)}$. 

% \section{Proof of Lemma \ref{lem:convex-max-relax}}\label{app:lem:convex-max-relax}
% When $\hat{Q}_k = Q_k$ is satisfied, there is
% \begin{equation*}
%     \frac{\alpha}{\|b\|_A} A + 2B =  Q_k\left[\frac{\alpha}{\|b\|_A}\hat{\Lambda}_k^{-2}\Lambda_k + \hat{\Lambda}_k^{-1}\right]Q_k^\top,
% \end{equation*}
% whose eigenvectors are the columns of $Q_k$. Lemma \ref{lem:convex-max-bstar} indicates that the optimal $b$ can be expressed as $\sqrt{\gamma_{0}(z_k)}\mathbf{v}_i$ for certain $i\in\{1,\ldots, d_x\}$, then we have the objective of \eqref{eq:convex-max} is $(2\sqrt{\gamma_{0}(z_k)\gamma_1\lambda_{i,k}}+\gamma_{0}(z_k))\hat{\lambda}_{i,k}^{-1}$. Since $j_k=\arg\max_{i}(2\sqrt{\gamma_{0}(z_k)\gamma_1\lambda_{i,k}}+\gamma_{0}(z_k))\hat{\lambda}_{i,k}^{-1}$, problem \eqref{eq:convex-max} is solved as
% \begin{equation*}
%     b^{*} = \sqrt{\gamma_{0}(z_k)}\mathbf{v}_{j_k}.
% \end{equation*}
\end{appendices}

\begin{IEEEbiography}[{\includegraphics[width=1in,height=1.5in,clip,keepaspectratio]{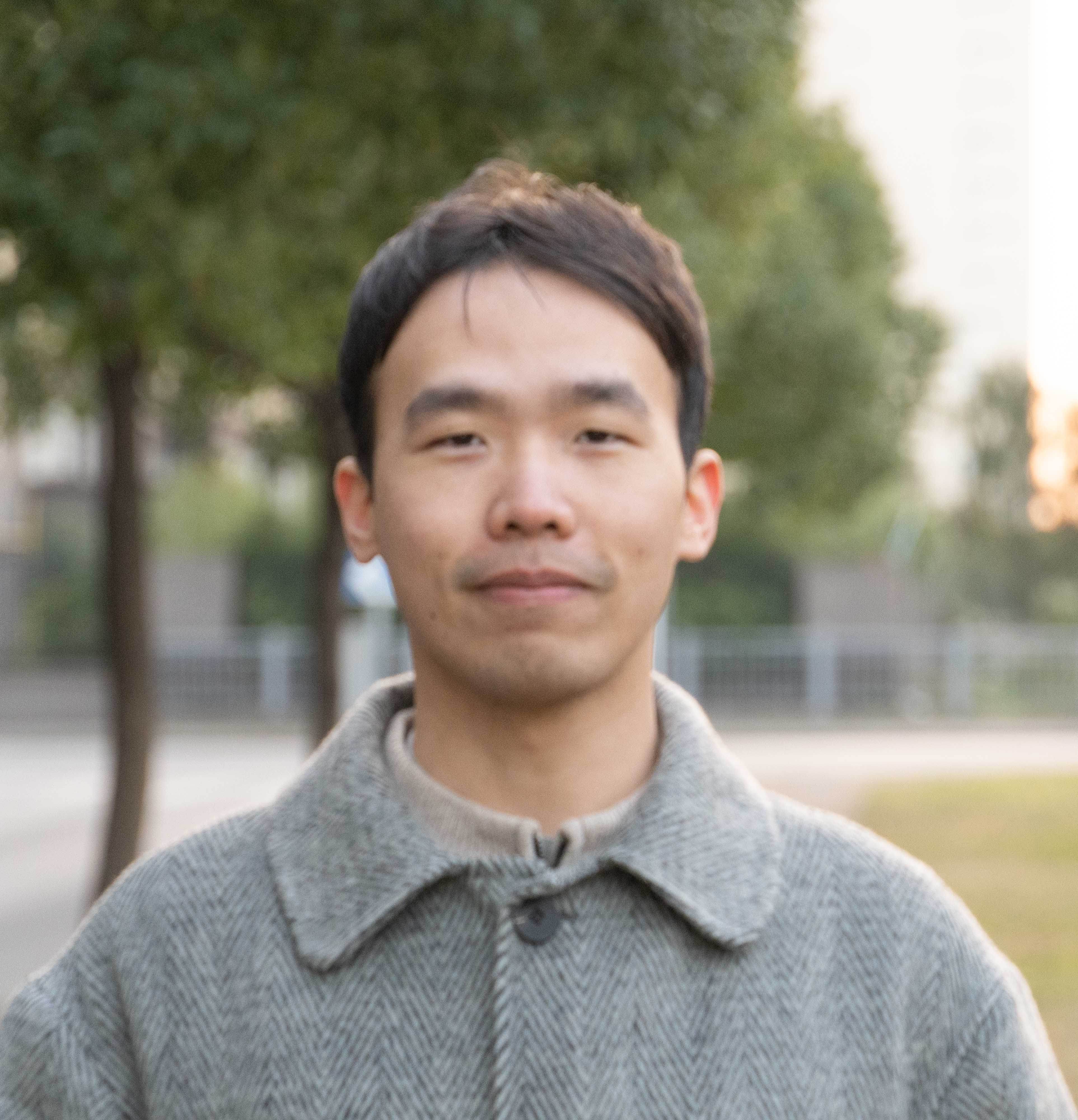}}]{Tao Xu}
	(S'22) received a B.S. degree in the School of Mathematical Sciences from Shanghai Jiao Tong University (SJTU), Shanghai, China. He is currently working toward the Ph.D. degree with the Department of Automation, SJTU. He is a member of Intelligent of Wireless Networking and Cooperative Control Group. His research interests include probabilistic prediction, distributionally robust optimization, dynamic games, and robotics.
\end{IEEEbiography}

\begin{IEEEbiography}[{\includegraphics[width=1in,height=1.5in,clip,keepaspectratio]{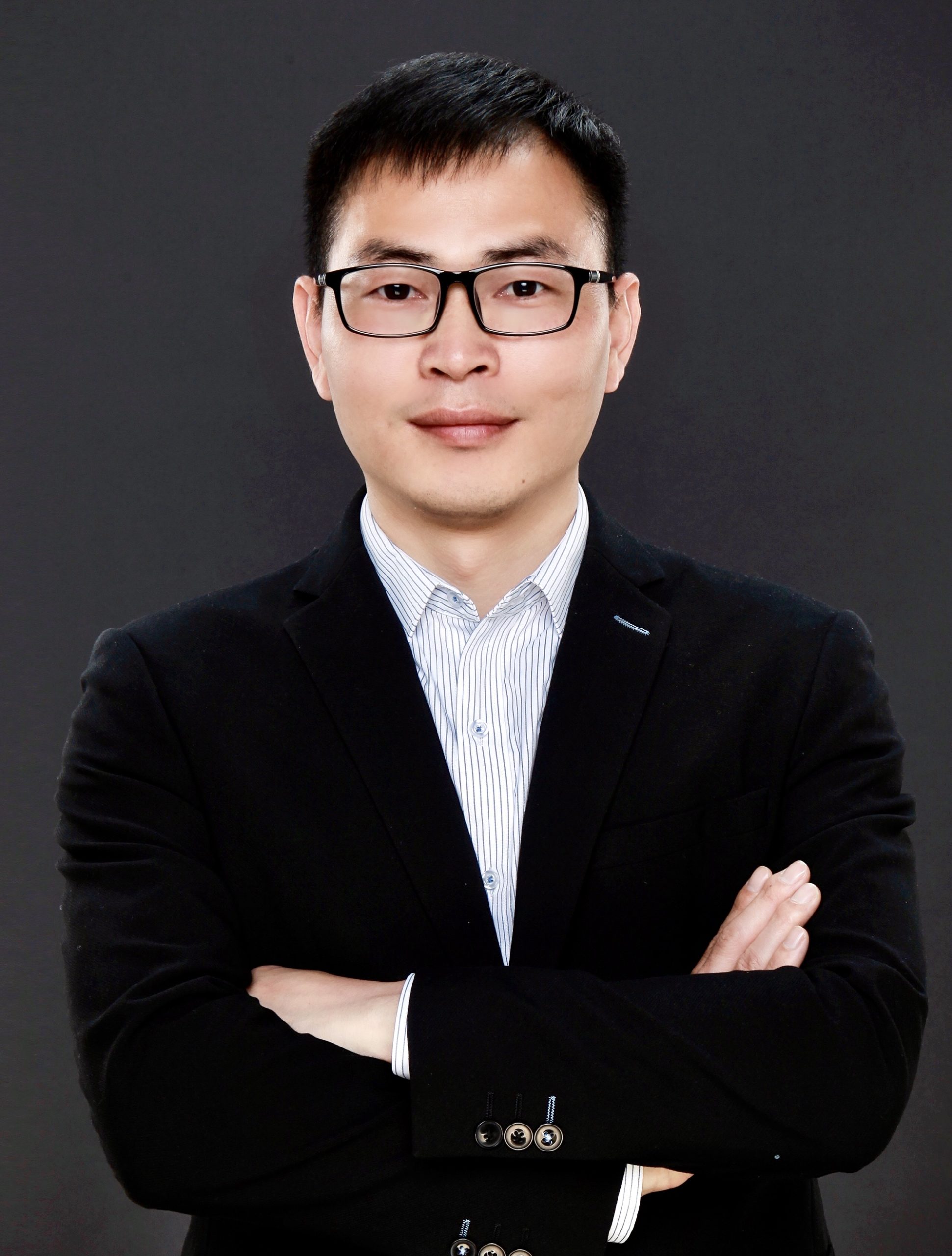}}]{Jianping He}
	(SM'19) is currently a full professor in the School of Automation and Intelligent Sensing at Shanghai Jiao Tong University, Shanghai, China. 
He received the Ph.D. degree in control science and engineering from Zhejiang University, Hangzhou, China, in 2013, and had been a research fellow in the Department of Electrical and Computer Engineering at University of Victoria, Canada, from Dec. 2013 to Mar. 2017. His research interests mainly include the distributed learning, control and optimization, security and privacy in network systems. 
\end{IEEEbiography}

\end{document}